\documentclass[12pt]{amsart}

\usepackage{amssymb}
\usepackage{amsmath}
\usepackage{amscd}
\usepackage{float}
\usepackage{graphicx}
\usepackage{amsfonts}
\usepackage{pb-diagram}
\usepackage{eufrak}

\newcommand{\Across}{\raisebox{-0.25\height}{\includegraphics[width=0.5cm]{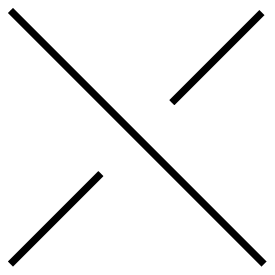}}}

\newcommand{\Asmooth}{\raisebox{-0.25\height}{\includegraphics[width=0.5cm]{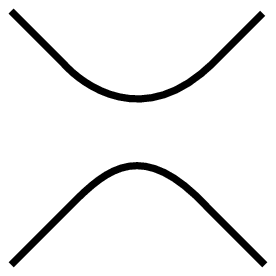}}}
\newcommand{\Bsmooth}{\raisebox{-0.25\height}{\includegraphics[width=0.5cm]{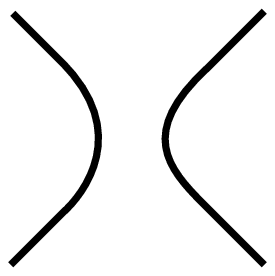}}}

\newcommand{\NegCross}{\raisebox{-0.25\height}{\includegraphics[width=0.5cm]{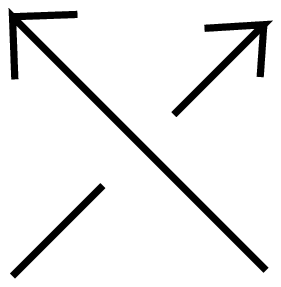}}}
\newcommand{\PosCross}{\raisebox{-0.25\height}{\includegraphics[width=0.5cm]{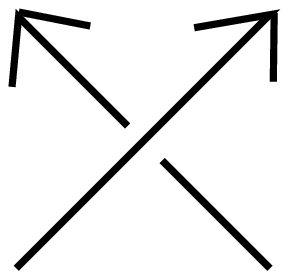}}}
\newcommand{\OrSmooth}{\raisebox{-0.25\height}{\includegraphics[width=0.5cm]{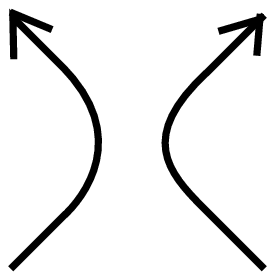}}}

\newcommand{\Rcurl}{\raisebox{-0.25\height}{\includegraphics[width=0.5cm]{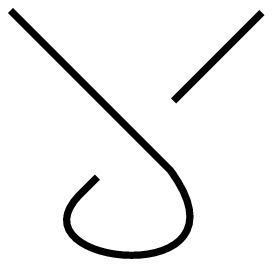}}}
\newcommand{\Lcurl}{\raisebox{-0.25\height}{\includegraphics[width=0.5cm]{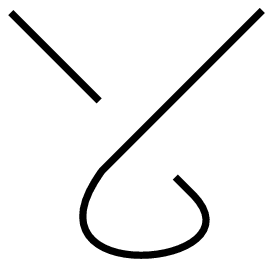}}}
\newcommand{\Arc}{\raisebox{-0.25\height}{\includegraphics[width=0.5cm]{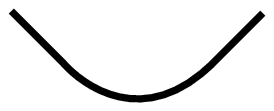}}}
\newcommand{\DegZero}{\raisebox{-0.25\height}{\includegraphics[width=1.0cm]{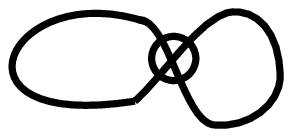}}}
\newcommand{\DegTwo}{\raisebox{-0.25\height}{\includegraphics[width=0.8cm]{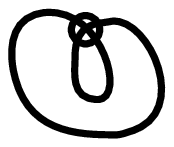}}}

\begin{document}

 \title[Rotational Virtual Knots and Quantum Link Invariants]{Rotational Virtual Knots and Quantum Link Invariants}

\author{Louis H. Kauffman}
\address{ Department of Mathematics, Statistics and
 Computer Science, University of Illinois at Chicago,
 851 South Morgan St., Chicago IL 60607-7045, U.S.A.}
\email{kauffman@math.uic.edu}
\urladdr{http://www.math.uic.edu/$~$kauffman/}



\date{}
\maketitle


\begin{abstract}
This paper studies rotational virtual knot theory and its relationship with quantum link invariants. Every quantum link invariant for classical knots and links extends to an invariant 
of rotational virtual knots and links. We give examples of non-trivial rotational virtual links that are undetectable by quantum invariants.
\end{abstract}

\section{Introduction}
This paper studies rotational virtual knot theory and its relationship with quantum link invariants. Every quantum link invariant for classical knots and links extends to an invariant 
of rotational virtual knots and links.\\

Virtual knot theory \cite{VKT,SVKT,DVK,VKCob} can be described diagrammatically in terms of classical diagrams with an extra crossing type (a virtual crossing) that is neither an over nor an under crossing.
The virtual crossing is an artifact of the planar or spherical representation of the virtual knot or link. Along with the Reidemeister moves, one adds a {\it detour move} that allows one to reroute 
consecutive sequences of virtual crossings. In standard virtual knot theory any rerouting is allowed. In {\it rotational virtual knot theory} the rerouting is restricted to regular homotopy of immersed curves
in the plane or on the surface of the two dimensional sphere. It turns out that this restriction is in exact accord with extending quantum link invariants for classical links to quantum link invariants for 
rotational virtual links. This means that rotational virtual links are a mathematical testing ground for the properties of quantum link invariants.\\

In this paper we illustrate a number of ideas and examples in the relationship of virtual, rotational virtual and quantum link invariants. Section 2 is a review of basics about virtual knot theory.
In Section 3 we define surface genus for both virtual links and rotational virtual links. We show how the Kauffman bracket polynomial naturally generalizes to a non-trivial rotational invariant both by combinatorial definition and in the context of its modeling as a quantum link invariant. We examine another bracket, the binary bracket and show how in this case it has a model as a quantum link invariant that in fact extends to an invariant of all virtual knots and links. Nevertheless, there is a combinatorial extension of the binary bracket to a non-trivial rotational invariant. We also show how the Manturov parity bracket for standard virtuals extends to an invariant of rotational virtuals, and we prove an {\it Irreducibilty Theorem} (Section 3.3) that says that if a virtual diagram has all of its crossings of odd parity, then it can be modified (by adding flat virtual curls) so that it is an irreducible non-trivial 
rotational diagram. This allows the calculation of the rotational genus of the modified diagram in terms of the diagram itself. Essentially, after modification, the odd diagrams become their own rotational
invariants. Then we go on to detail the construction of oriented quantum link invariants for both classical and rotational virtual knots and links (Section 5) and end this section with the extension of 
the Homflypt polynomial. A question is raised about the relationship of these extensions with the Khovanov-Rozansky categorification of specializations of the Homflypt polynomial.\\

The paper ends with a section on quantum link invariants in the Hopf algebra framework where one can see the naturality of using regular homotopy combined with virtual crossings (permutation operators), as they occur significantly in the category associated with a Hopf algebra. We show how this approach via categories and quantum algebras illuminates the structure of invariants that we have already described via state summations. In particular, we show in Figure~\ref{nontriv} how the link of Figure~\ref{linkcalc2} gets a non-trivial functorial image, corresponding to its non-trivial rotational bracket invariant. We also show that the link of Figure~\ref{linkcalc3} has, in Figure~\ref{triv},  a trivial functorial image. This means that not only is this link not detected by the rotational bracket polynomial, it is not detected by any quantum invariant formulated as outlined in this section. The link of Figure~\ref{linkcalc3} is a non-trivial rotational link, as we show in Section 3.3 by using a parity technique. These calculations with the diagrammatic images in quantum algebra category show how this category forms a higher level language for understanding rotational virtual knots and links. There are inherent limitations to studying rotational virtual knots by quantum algebra alone.
\\

This paper is partly expository in its intent and contains many new results. Our intention is that this paper can be a springboard for many further investigations into rotational virtual knot theory
and quantum invariants.\\

 \section{Virtual Knot Theory}
Knot theory
studies the embeddings of curves in three-dimensional space.  Virtual knot theory studies the  embeddings of curves in thickened surfaces of arbitrary
genus, up to the addition and removal of empty handles from the surface. Virtual knots have a special diagrammatic theory, described below. Classical knot
theory embeds in virtual knot theory.
\bigbreak  

\begin{figure}[htb]
     \begin{center}
     \begin{tabular}{c}
     \includegraphics[width=8cm]{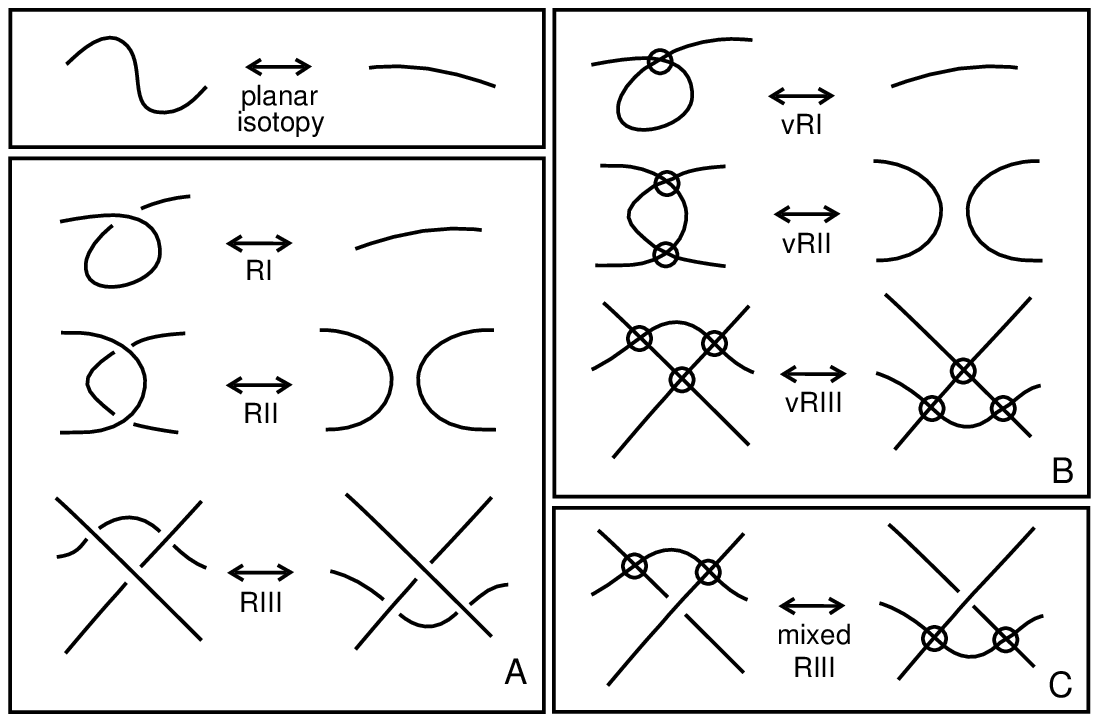}
     \end{tabular}
     \caption{\bf Moves}
     \label{Figure 1}
\end{center}
\end{figure}

In the diagrammatic theory of virtual knots one adds 
a {\em virtual crossing} (see Figure~\ref{Figure 1}) that is neither an over-crossing
nor an under-crossing.  A virtual crossing is represented by two crossing segments with a small circle
placed around the crossing point. 
\bigbreak

\begin{figure}[htb]
     \begin{center}
     \begin{tabular}{c}
     \includegraphics[width=6cm]{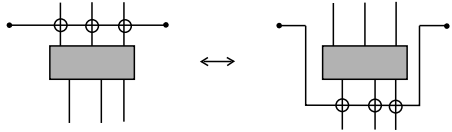}
     \end{tabular}
     \caption{\bf Detour Move}
     \label{Figure 2}
\end{center}
\end{figure}

Moves on virtual diagrams generalize the Reidemeister moves for classical knot and link
diagrams.  See  Figure~\ref{Figure 1}.  Classical crossings interact with
one another according to the usual Reidemeister moves, while virtual crossings are artifacts of the structure in the plane. 
Adding the global detour move to the Reidemeister moves completes the description of moves on virtual diagrams. In  Figure~\ref{Figure 1} we illustrate a set of local
moves involving virtual crossings. The global detour move is
a consequence of  moves (B) and (C) in  Figure~\ref{Figure 1}. The detour move is illustrated in  Figure~\ref{Figure 2}.  Virtual knot and link diagrams that can be connected by a finite 
sequence of these moves are said to be {\it equivalent} or {\it virtually isotopic}.
\bigbreak

\begin{figure}[htb]
     \begin{center}
     \begin{tabular}{c}
     \includegraphics[width=6cm]{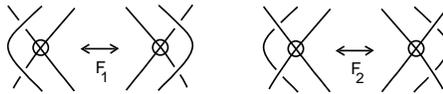}
     \end{tabular}
     \caption{\bf Forbidden Moves}
     \label{Figure 3}
\end{center}
\end{figure}

Another way to understand virtual diagrams is to regard them as representatives for oriented Gauss codes \cite{GPV}, \cite{VKT,SVKT} 
(Gauss diagrams). Such codes do not always have planar realizations. An attempt to embed such a code in the plane
leads to the production of the virtual crossings. The detour move makes the particular choice of virtual crossings 
irrelevant. {\it Virtual isotopy is the same as the equivalence relation generated on the collection
of oriented Gauss codes by abstract Reidemeister moves on these codes.}  
\bigbreak

 Figure~\ref{Figure 3} illustrates the two {\it forbidden moves}. Neither of these follows from Reidmeister moves plus detour move, and 
indeed it is not hard to construct examples of virtual knots that are non-trivial, but will become unknotted on the application of 
one or both of the forbidden moves. The forbidden moves change the structure of the Gauss code and, if desired, must be 
considered separately from the virtual knot theory proper. 
\bigbreak

\subsection{Interpretation of Virtuals Links as Stable Classes of Links in  Thickened Surfaces}
\begin{figure}[htb]
     \begin{center}
     \begin{tabular}{c}
     \includegraphics[width=8cm]{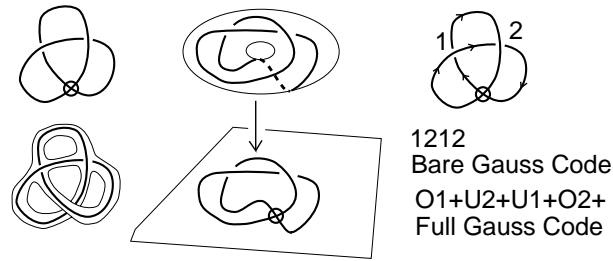}
     \end{tabular}
     \caption{\bf Surface Representation}
     \label{Figure 4}
\end{center}
\end{figure}

There is a useful topological interpretation \cite{VKT,DVK,CS1} for this virtual theory in terms of embeddings of links
in thickened surfaces.  Regard each 
virtual crossing as a shorthand for a detour of one of the arcs in the crossing through a 1-handle
that has been attached to the 2-sphere of the original diagram.  
By interpreting each virtual crossing in this way, we
obtain an embedding of a collection of circles into a thickened surface  $S_{g} \times R$ where $g$ is the 
number of virtual crossings in the original diagram $L$, $S_{g}$ is a compact oriented surface of genus $g$
and $R$ denotes the real line.  We say that two such surface embeddings are
{\em stably equivalent} if one can be obtained from another by isotopy in the thickened surfaces, 
homeomorphisms of the surfaces and the addition or subtraction of empty handles (i.e. the knot does not go through the handle).
See Figure~\ref{Figure 4} for an illustration of the relationship of virtual diagrams with diagrams in surfaces, and with the Gauss code representation.\\

\noindent We have the
\smallbreak
\noindent
{\bf Theorem 1 \cite{VKT,DVK,DVK,CS1}.} {\em Two virtual link diagrams are isotopic if and only if their corresponding 
surface embeddings are stably equivalent.}  
\smallbreak
\noindent
\bigbreak  

\noindent  The reader will find more information about this correspondence \cite{VKT,DVK} in other papers by the author and in the literature of virtual knot theory.
\bigbreak
 
  \subsection{Review of the Bracket Polynomial for Virtual Knots}

In this section we recall how the bracket state summation model \cite{KaB} for the Jones polynomial  is defined for virtual knots
and links.  
\bigbreak

The bracket polynomial \cite{KaB} model for the Jones polynomial \cite{JO,JO1,JO2,Witten} is usually described by the expansion
\begin{equation}
\langle \Across \rangle=A \langle \Asmooth \rangle + A^{-1}\langle
\Bsmooth \rangle \label{kabr}
\end{equation}

and we have

\begin{equation}
\langle K \, \bigcirc \rangle=(-A^{2} -A^{-2}) \langle K \rangle \label{kabr1}
\end{equation}

\begin{equation}
\langle \Rcurl \rangle=(-A^{3}) \langle \Arc \rangle \label{kabr2}
\end{equation}

\begin{equation}
\langle \Lcurl \rangle=(-A^{-3}) \langle \Arc \rangle \label{kabr3}
\end{equation}
\bigbreak

We call a diagram in the plane 
{\em purely virtual} if the only crossings in the diagram are virtual crossings. Each purely virtual diagram is equivalent by the
virtual moves to a disjoint collection of circles in the plane.
\bigbreak

A state $S$ of a link diagram $K$ is obtained by
choosing a smoothing for each crossing in the diagram and labelling that smoothing with either $A$ or $A^{-1}$
according to the convention indicated in the bracket expansion above.  Then, given
a state $S$, one has the evaluation $<K|S>$ equal to the product of the labels at the smoothings, and one has the 
evaluation $||S||$ equal to the number of loops in the state (the smoothings produce purely virtual diagrams).  One then has
the formula
$$<K> = \Sigma_{S}<K|S>d^{||S||-1}$$
where the summation runs over the states $S$ of the diagram $K$, and $d = -A^{2} - A^{-2}.$
This state summation is invariant under all classical and virtual moves except the first Reidemeister move.
The bracket polynomial is normalized to an
invariant $f_{K}(A)$ of all the moves by the formula  $f_{K}(A) = (-A^{3})^{-w(K)}<K>$ where $w(K)$ is the
writhe of the (now) oriented diagram $K$. The writhe is the sum of the orientation signs ($\pm 1)$ of the 
crossings of the diagram. The Jones polynomial, $V_{K}(t)$ is given in terms of this model by the formula
$$V_{K}(t) = f_{K}(t^{-1/4}).$$
\noindent This definition is a direct generalization to the virtual category of the  
state sum model for the original Jones polynomial. It is straightforward to verify the invariances stated above.
In this way one has the Jones polynomial for virtual knots and links.
\bigbreak

\noindent We have \cite{DVK} the  
\smallbreak
\noindent
{\bf Theorem.} {\em To each non-trivial
classical knot diagram of one component $K$ there is a corresponding  non-trivial virtual knot diagram $Virt(K)$ with unit
Jones polynomial.} 
\bigbreak

\begin{figure}[htb]
     \begin{center}
     \begin{tabular}{c}
     \includegraphics[width=7cm]{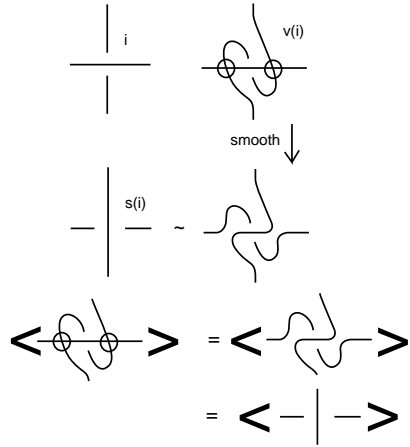}
     \end{tabular}
     \caption{\bf Virtualizing a Crossing and Crossing Switches}
     \label{Figure 5}
\end{center}
\end{figure}

\begin{figure}[htb]
     \begin{center}
     \begin{tabular}{c}
     \includegraphics[width=5cm]{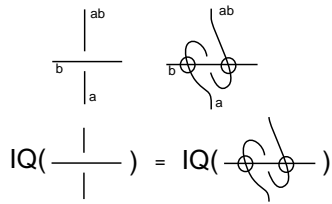}
     \end{tabular}
     \caption{\bf Quandle Invariance Under Virtualization}
     \label{Figure 6}
\end{center}
\end{figure}

The main ideas behind this theorem are indicated in  Figure~\ref{Figure 5} and  Figure~\ref{Figure 6}.  In  Figure~\ref{Figure 5} we indicate the 
virtualization operation that replaces a classical crossing by using two virtual crossings and changing the implicit orientation of the classical crossing. We also show how the bracket polynomial sees this operation as though the crossing had been switched in the classical knot. Thus, if we virtualize the set
of classical crossings whose switching will unknot the knot, then the virtualized knot will have unit 
Jones polynomial. On the other hand, the virtualization is invisible to the quandle, as shown in  
Figure~\ref{Figure 6}. This implies (by properties of the quandle) that virtual knots obtained in this way from classical non-trivial knots will themselves be non-trivial.
In fact, such knots are non-classical \cite{DKK}.
\bigbreak

It is an open problem whether there are classical knots (actually knotted) having unit Jones polynomial.  Examples produced by virtualization are guaranteed to be non-trivial. It is 
more difficult to prove in general \cite{DKK} or even in special cases that such virtualizations are non-classical.  This has led to the discovery of  new invariants for virtual knots.\\

\section{Rotational Virtual Knot Theory}
In {\it rotational virtual knot theory} (introduced in~\cite{VKT} ) the detour move is restricted to regular homotopy of plane curves. This means that the virtual curl of Figure~\ref{Figure 1} can not be directly simplified, but two opposite virtual curls can be created or destroyed by using the Whitney Trick of Figure~\ref{whitney}. Another way to put this is to say that rotational virtual knot theory is virtual
knot theory without the first virtual move (thus one does not allow the addition or deletion of a virtual curl).\\

If we wish to take advantage of the integral Whitney degree of plane curves, then diagrams are represented in the surface of plane so that we can distinguish clockwise from counterclockwise rotations.
For many purposes this is preferable. In some cases we use the $2$--sphere, and use the equivalence classes of immersed curves up to regular isotopy in the $2$--sphere. These classes are in correspondence with the absolute values of the Whitney degree of corresponding plane curve immersions. \\

See Figure~\ref{rotbracket} for an
example of a rotational virtual knot.
The rotational  version of virtual knot theory is significant because {\it  all quantum link invariants
originally defined for classical links extend to rotational virtual
knot theory.}  We sketch this relationship now and then devote the later sections of the paper to it in more detail. See Figure~\ref{fig20} where we indicate how a quantum link invariant depends on matrices or operators assigned to each crossing and each maximum and minimum of the diagram. One extends this to virtual crossings by using a crossed identity operator ( a transposition) at the virtual crossings. A maximum or minimum composed with a permutation (virtual crossing) forms a virtual curl. Since a non-trivial operator for a maximum or a minimum may not be invariant under permutation, the quantum calculation may not be invariant under the removal of virtual curls. This theory has been begun in \cite{VKT,CVBraid}. We carry it further here.  In the next section we discuss quantum link invariants from the beginning.\\

The extension by a crossed identity operator for virtual crossings is the simplest way to handle a quantum link invariant for virtual knots and links. We can use more complex operators at a virtual crossing, and sometimes the addition of such a virtual operator can transform the invariant into one that is invariant under the first virtual move. We will discuss examples of this kind at the end of  Section 5. Thus in some cases we obtain invariants of virtual knots by extension from quantum invariants of classical knots.\\

\subsection{Genus of Rotational Virtual Knots and Links}
\begin{figure}[htb]
     \begin{center}
     \begin{tabular}{c}
     \includegraphics[width=8cm]{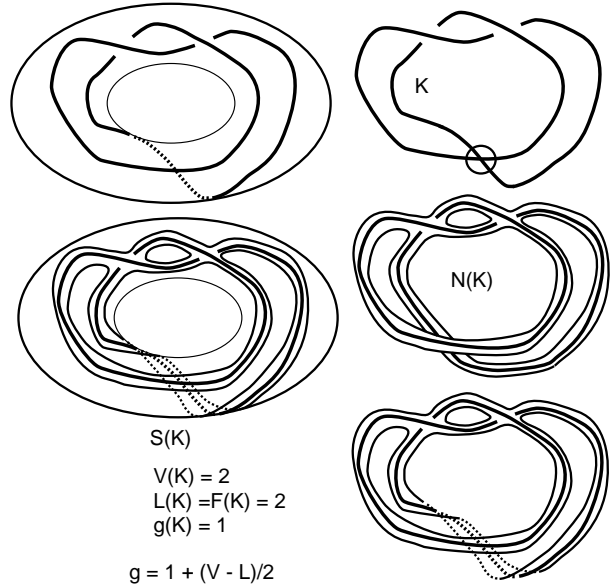}
     \end{tabular}
     \caption{\bf Genus of Virtual Diagrams}
     \label{genus}
\end{center}
\end{figure}

Since standard virtual knots and links can be represented as classical diagrams on oriented surfaces (or equivalently as embeddings in thickened surfaces) every virtual knot or link has a representation 
of least genus. Classical knots and links are those with least genus equal to zero. In Figure~\ref{genus} we illustrate a process for estimating the genus of a virtual link, by determining the genus of its
diagram. In this figure we show how to form a ribbon neighborhood $N(K)$ of a given virtual diagram $K$ that makes $K$ a classical diagram in the surface formed by $N(K).$ This surface is a surface with $L(K)$ boundary components. We can then add discs to the boundary components of $N(K)$ to form a closed surface $S(K)$ in which $K$ lives as a classical diagram. We {\it define} for the {\it diagram} $K$, the {\it diagram genus} $g(K) = genus(S(K)).$ The {\it virtual genus} $vg(K)$ is the minimum of $g(D)$ over all (virtual) diagrams $D$ that are virtually isotopic to $K.$ It is of interest to 
note that $$g(K) = 1 + (V(K) - L(K))/2,$$ where $V(K)$ denotes the number of classical crossings in the diagram $K$ and $L(K)$ is the number of boundary loops in $N(K).$ We omit the proof of this result. It is an Euler characteristic calculation.\\

\begin{figure}[htb]
     \begin{center}
     \begin{tabular}{c}
     \includegraphics[width=8cm]{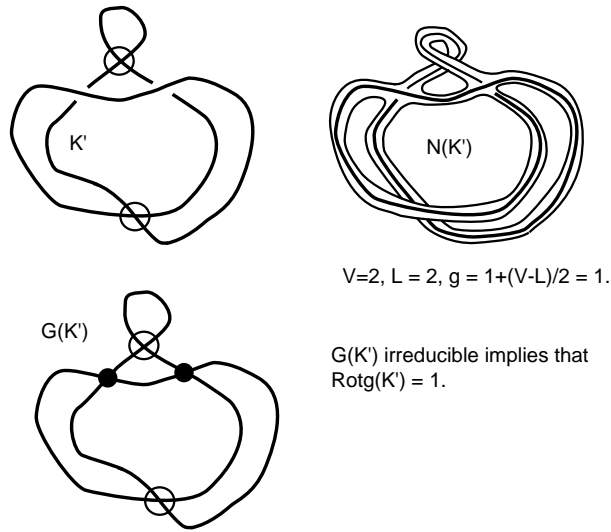}
     \end{tabular}
     \caption{\bf Genus of Rotational Virtual Diagrams}
     \label{genusrot}
\end{center}
\end{figure}

If $K$ is regarded as a rotational virtual link diagram, then $g(K) = 1 + (V(K) - L(K))/2$ is still defined by the same process of forming $N(K)$ as in Figure~\ref{genus} , but the set of diagrams $D$ that 
are rotationally equivalent to $K$ is putatively smaller than the full set of diagrams virtually isotopic to $K.$ Accordingly, we define the {\it rotational genus} of a virtual diagram $K$ to be 
$Rg(K),$ the minimum of $g(D)$ over all virtual diagrams $D$ that are rotationally equivalent to $K.$ In Figure~\ref{genusrot} we illustrate a virtual diagram $K'$ that is rotationally non-trivial, as we shall prove in the next section and has $Rg(K') =1.$ As a virtual knot $K'$ is trivial, but as a rotational virtual knot, $K'$ is non-trivial and has genus equal to one. Later in the paper, we will see other examples where one can determine the rotational genus of a virtual diagram.\\

\noindent {\bf Remark.} It is worth noting that when we convert a rotational virtual diagram to a surface with boundary and embedded link (a so-called {\it abstract link diagram}) as in Figure~\ref{genusrot} then 
the bands in the abstract link diagram may acquire curls as indicated in that Figure. If we keep track of these curls, as a kind of framing on the abstract link diagrams, then this is equivalent to working in the 
rotational virtual category. Furthermore, as in Figure~\ref{genus}, an abstract link diagram obtained from a particular embedded surface representation of a virtual knot may acquire a similar framing structure.
This points to geometric and framing possibilities for rotational virtuals. This approach will be taken up in a separate publication.\\

\subsection{The Rotational Bracket Polynomial}
In this section we shall examine a direct model of the  the bracket polynomial that gives information about rotational virtual knots and links.
We  formulate \cite{VKT} a version of the bracket polynomial for rotational
virtuals that assigns variables according to the Whitney degree of state curves. 
For rotational virtuals we extend the bracket just as we did for virtual knots and links except that the state curves are now disjoint unions in the plane of curves that have only virtual self-intersections, and are taken up to {\it regular homotopy} in the $2$--sphere (we can also use the plane, but for convenience we use the $2$--sphere in this section and so can use unoriented diagrams). We can take regular homotopy of curves to mean that the flat virtual versions of the second and third Reidemeister moves are allowed; regular homotopy is the equivalence relation generated by these moves and planar homeomorphisms.\\

\begin{figure}
     \begin{center}
     \begin{tabular}{c}
     \includegraphics[width=8cm]{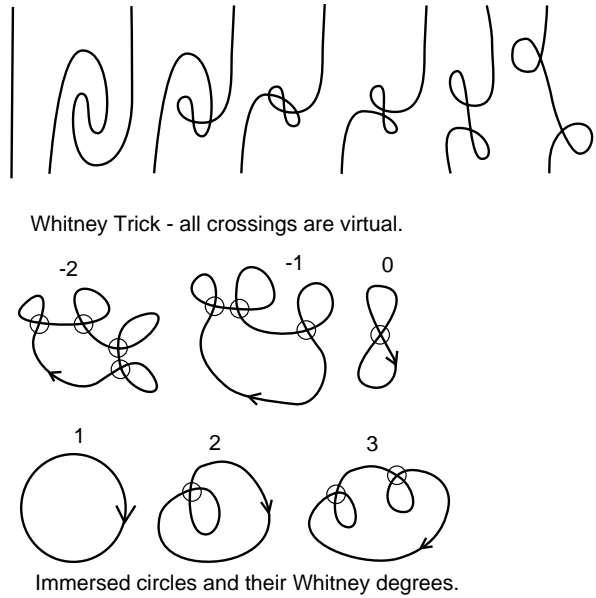}
     \end{tabular}
     \caption{\bf Whitney Trick and Whitney Degrees}
     \label{whitney}
\end{center}
\end{figure}

We shall use the same notation, $<K>$, for the rotational bracket, but it is understood that the detour moves are taken up to regular isotopy and that there is no given invariance under removal or addition of a virtual curl. When we expand the bracket we obtain a state sum of the form
$$<K> = \Sigma_{S}<K|S>[S]$$
where the summation is over all states obtained by smoothing every crossing in the virtual diagram $K$
and $<K|S>$ is the product of the weights $A$ and $A^{-1}$ just as before. An empty loop with no virtual crossings (in its virtual equivalence class) will be evaluated as $d = -A^2 - A^{-2}.$
The symbol  $[S]$ is the {\it class} of the state $S.$ By the class of the state we mean its equivalence class up to virtual rotational equivalence. This means that each state loop is taken as a regular homotopy class (in the $2$--sphere). These individual classes are in $1$-$1$ correspondence with the non-negative integers (for the $2$--sphere), as shown in Figure~\ref{whitney} (via the Whitney trick and the winding degree of the plane curves), and can be handled by using combinatorial regular isotopy as in \cite{FKT}.  A configuration of loops (possibly nested) is equivalent to a disjoint union of adjacent loops. We can thus regard each virtual loop as a variable $d_{n}$ where $n$ is $0$ or a positive integer and $d_{1} =d_{-1} = -A^2 - A^{-2}.$ We may use diagrams for $d_{0}$ and $d_{2}$ as shown below:
$$d_{0} = \DegZero, d_{2} = \DegTwo.$$

\begin{figure}
     \begin{center}
     \begin{tabular}{c}
     \includegraphics[width=6cm]{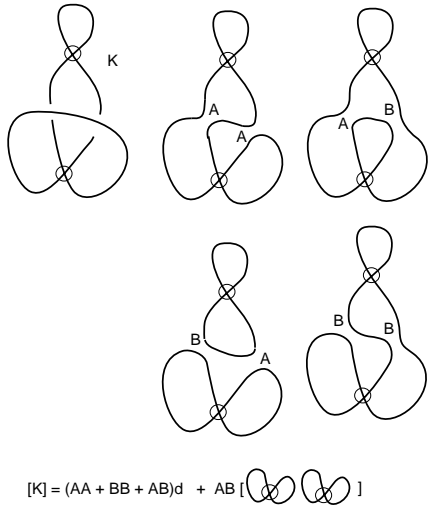}
     \end{tabular}
     \caption{\bf Bracket Expansion of a Rotational Virtual Knot}
     \label{rotbracket}
\end{center}
\end{figure}

Here we give examples of a computation of 
$<K>$ for a rotational virtual knot in Figure~\ref{rotbracket}. For the knot $K$ in this figure, we show the computation of the raw rotational bracket with variables $A,B,d.$ The result is
$$<K> = (A^2 + B^2 + AB)d + AB[\DegZero]^2 =  (A^2 + B^2 + AB)d + AB d_{0}^{2} .$$
The reader will note that in this example, even if we let $A = -1 = B$ and $d = -2$ the invariant is still non-trivial due to the appearance of the two loops with Whitney degree zero. Thus the example in Figure~\ref{rotbracket} also gives a non-trivial flat rotational virtual knot. In \cite{VKCob} we consider cobordism of rotational virtual knots, but we shall not discuss that aspect in this paper.
\bigbreak

\begin{figure}
     \begin{center}
     \begin{tabular}{c}
     \includegraphics[width=6cm]{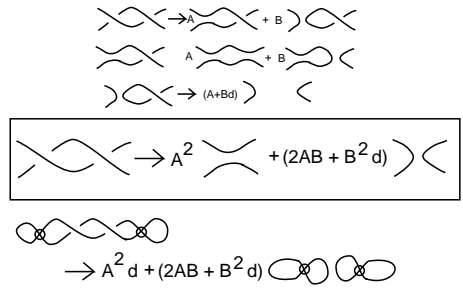}
     \end{tabular}
     \caption{\bf Expanding a Tangle}
     \label{tanglexp}
\end{center}
\end{figure}

\begin{figure}
     \begin{center}
     \begin{tabular}{c}
     \includegraphics[width=6cm]{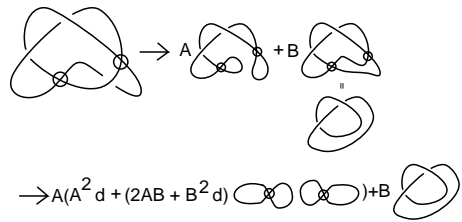}
     \end{tabular}
     \caption{\bf Evaluation of Rotational Bracket on a Virtual Knot}
     \label{knotexp}
\end{center}
\end{figure}

In Figures ~\ref{tanglexp} and ~\ref{knotexp} we show the calculation of the raw rotational bracket for a virtual knot whose ordinary normalized bracket is one. (It is a virtual knot with unit Jones polynomial.)
We see that the result gives a non trivial rotational bracket when we take $B = A^{-1}$ and $d = -A^2 - A^{-2}.$ \\

\begin{figure}
     \begin{center}
     \begin{tabular}{c}
     \includegraphics[width=6cm]{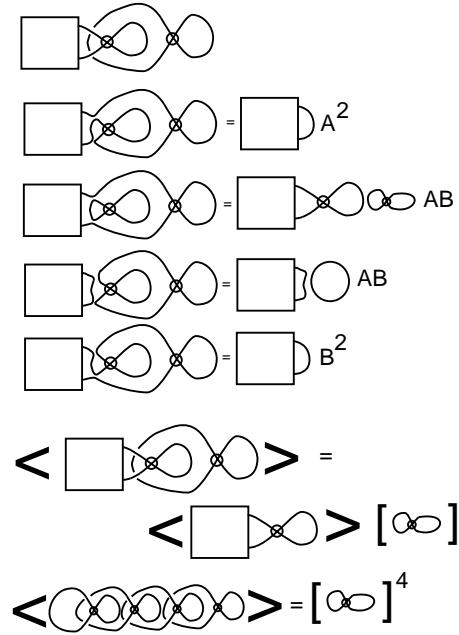}
     \end{tabular}
     \caption{\bf Evaluation of Rotational Bracket on a Family of Links}
     \label{linkrecur}
\end{center}
\end{figure}

In Figure~\ref{linkrecur} we show how to calculate the rotational bracket for a family of links and we show that for the choice  $B = A^{-1}$ and $d = -A^2 - A^{-2}$ the polynomials are of the form
$$<L_{n}> = [\DegZero]^{n}.$$ This means that the rotational bracket does not distinguish these links from unlinks with the same components. We shall return to the problem of detecting these links later in the paper.\\

\begin{figure}
     \begin{center}
     \begin{tabular}{c}
     \includegraphics[width=6cm]{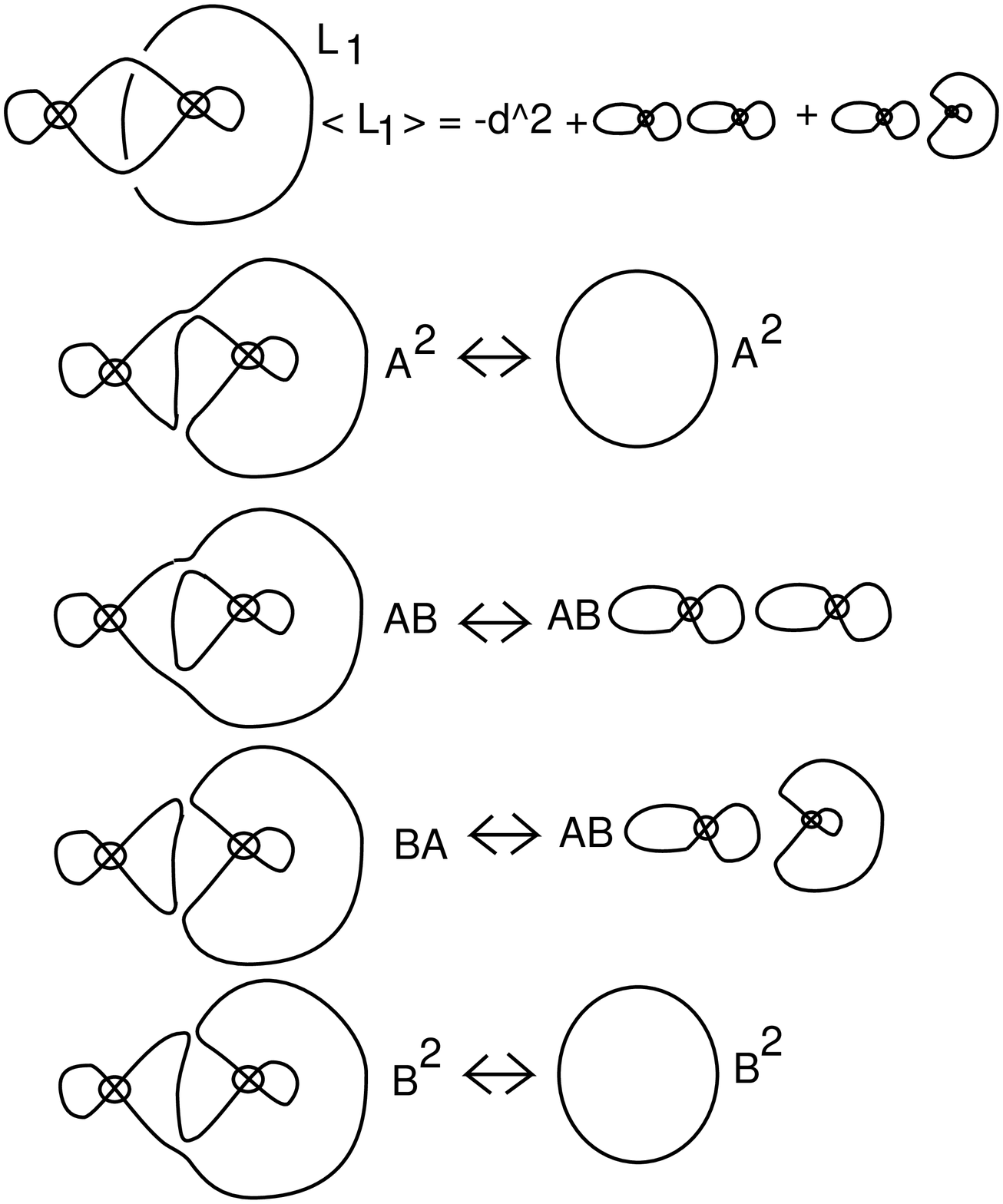}
     \end{tabular}
     \caption{\bf Evaluation of Rotational Bracket on a Link $L_1$}
     \label{linkcalc1}
\end{center}
\end{figure}

\begin{figure}
     \begin{center}
     \begin{tabular}{c}
     \includegraphics[width=6cm]{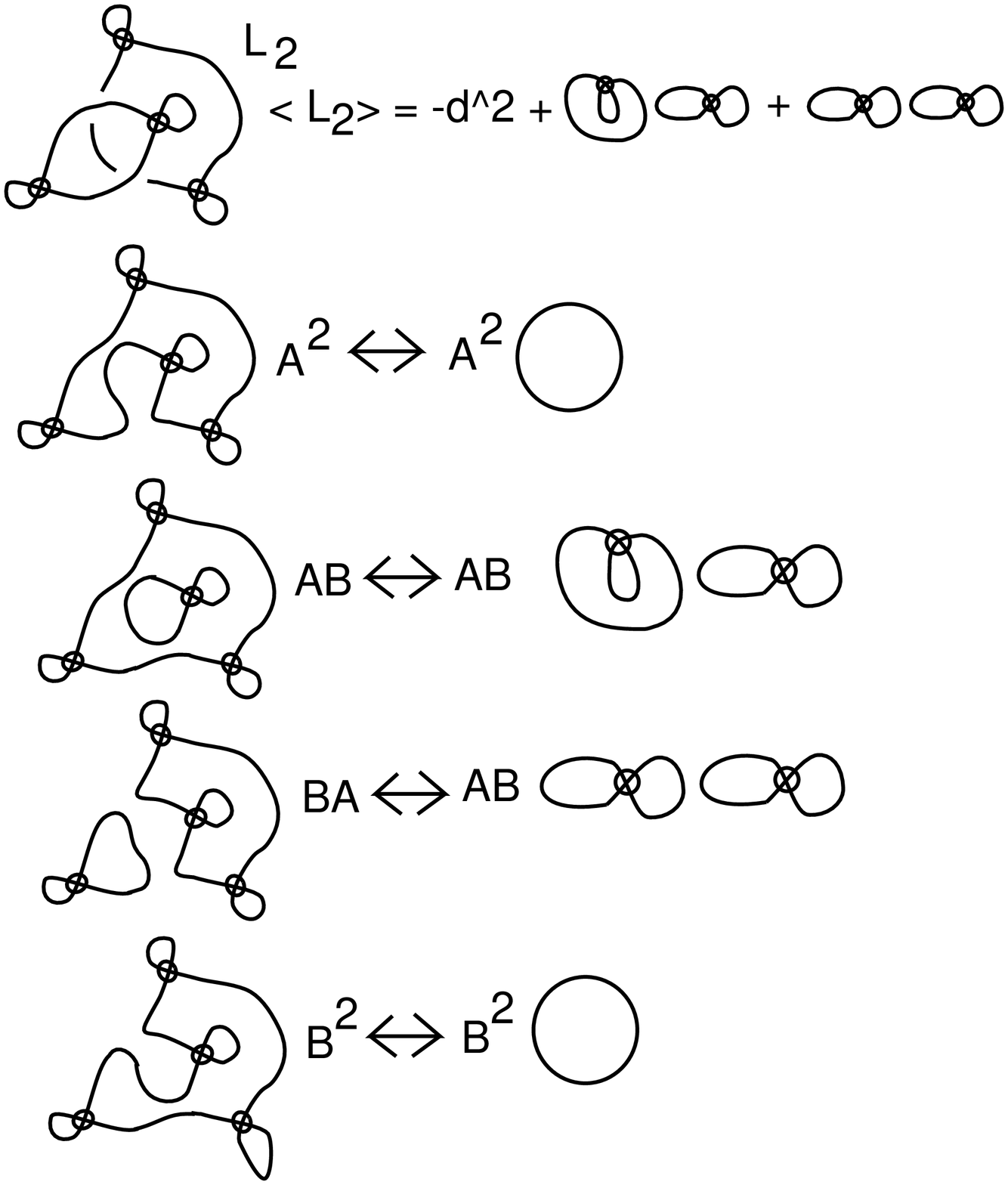}
     \end{tabular}
     \caption{\bf Evaluation of Rotational Bracket on a Link $L_2$}
     \label{linkcalc2}
\end{center}
\end{figure}

\begin{figure}
     \begin{center}
     \begin{tabular}{c}
     \includegraphics[width=6cm]{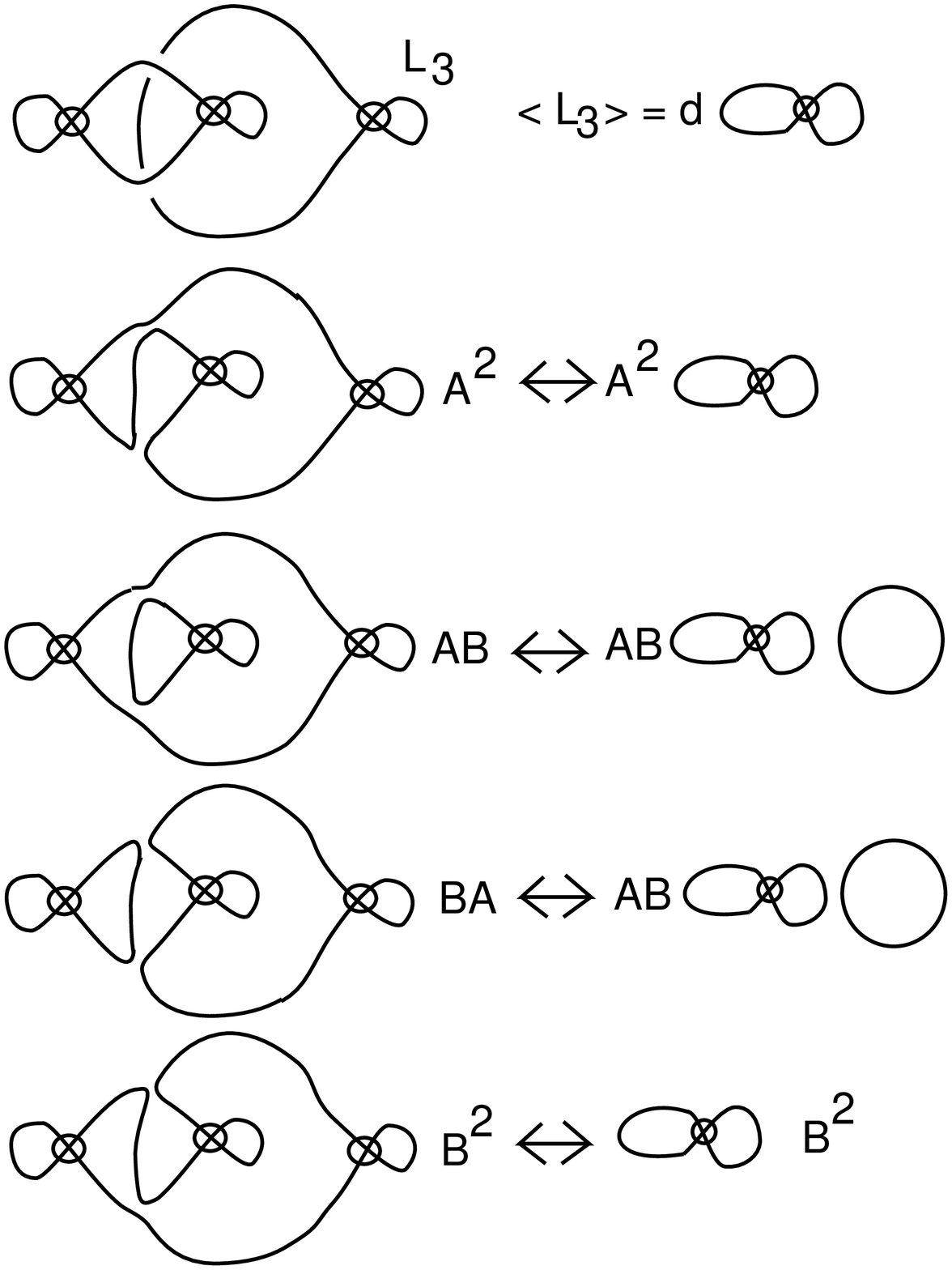}
     \end{tabular}
     \caption{\bf Evaluation of Rotational Bracket on a Link $L_3$}
     \label{linkcalc3}
\end{center}
\end{figure}

\begin{figure}
     \begin{center}
     \begin{tabular}{c}
     \includegraphics[width=6cm]{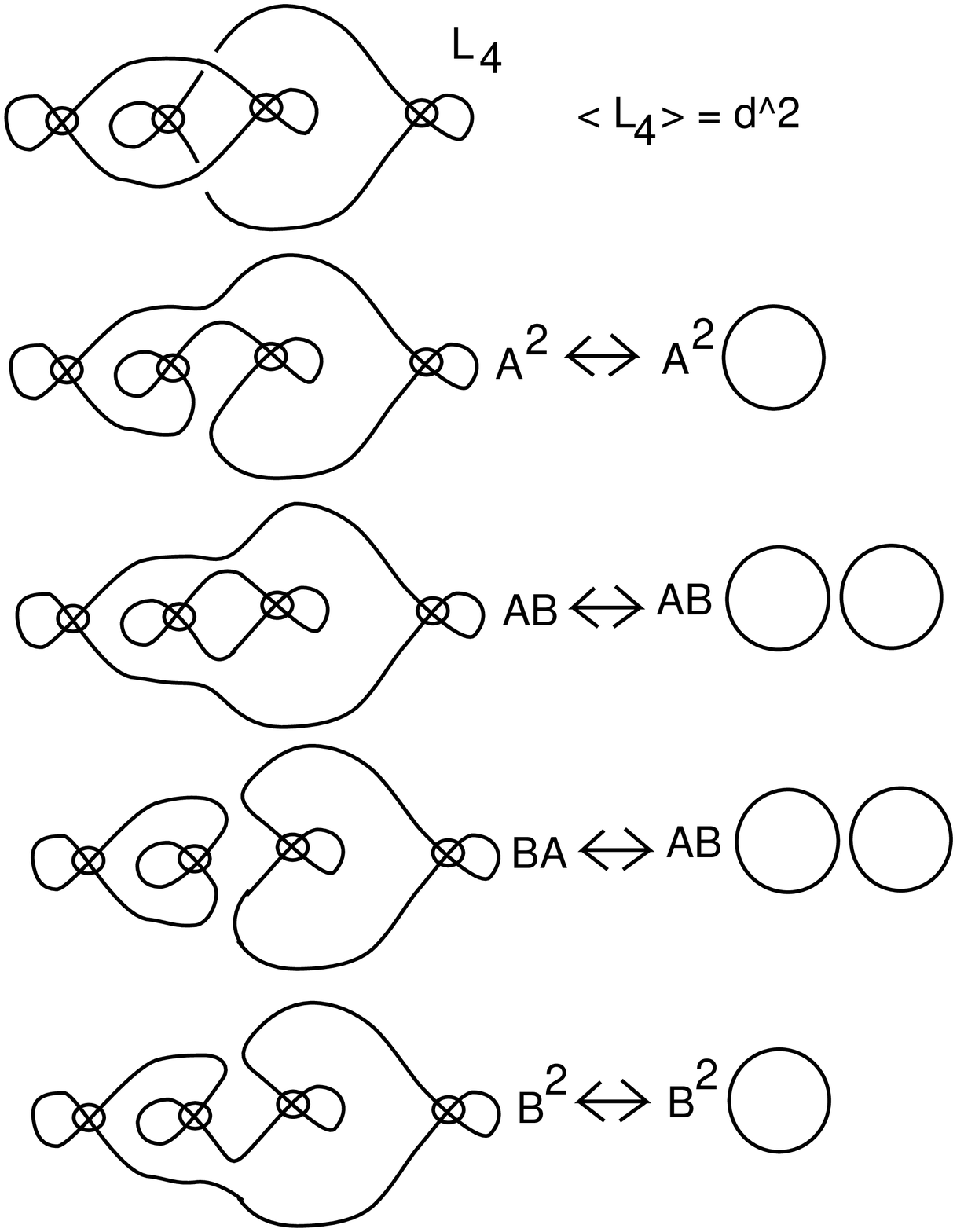}
     \end{tabular}
     \caption{\bf Evaluation of Rotational Bracket on a Link $L_4$}
     \label{linkcalc4}
\end{center}
\end{figure}

\begin{figure}
     \begin{center}
     \begin{tabular}{c}
     \includegraphics[width=6cm]{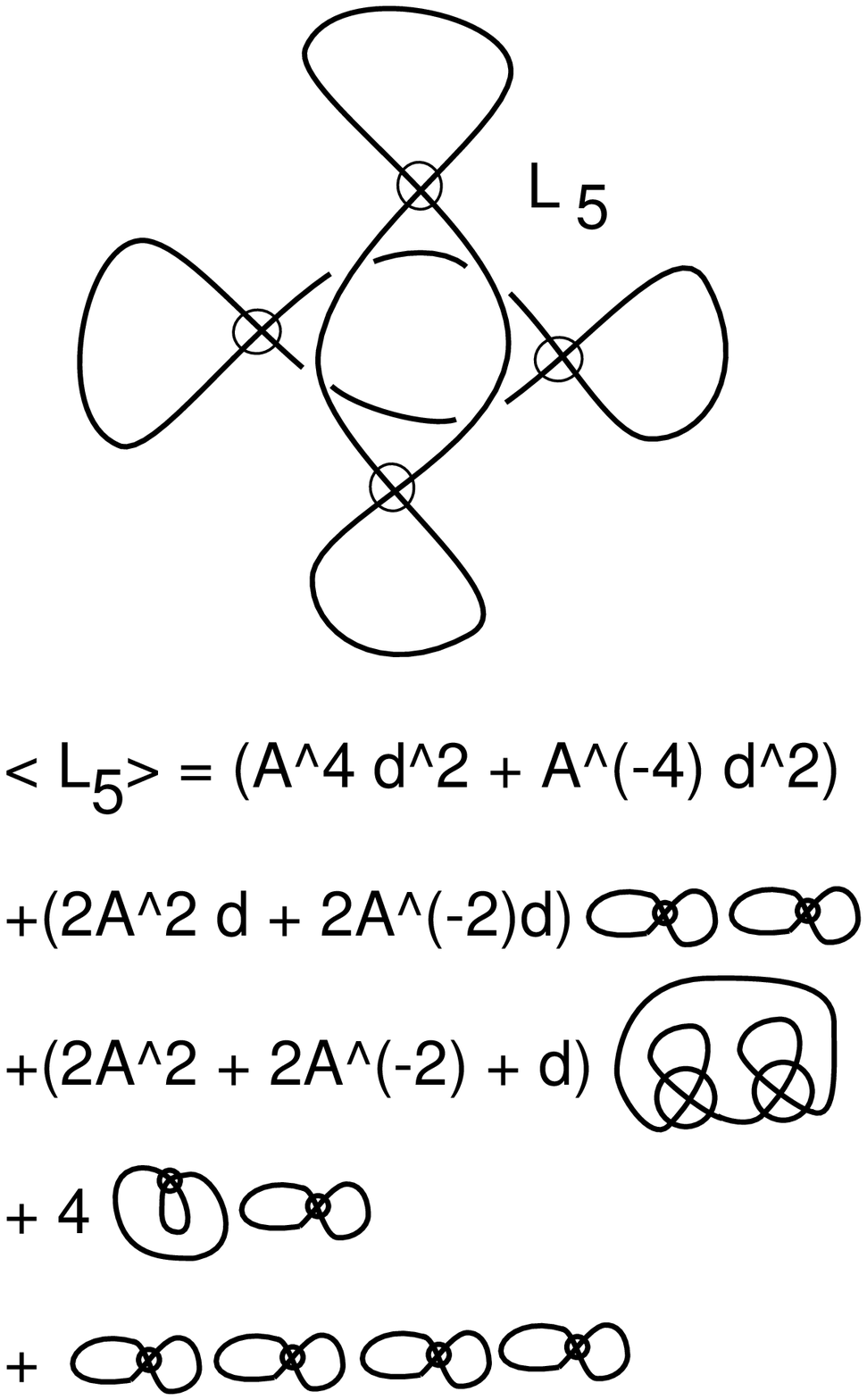}
     \end{tabular}
     \caption{\bf Evaluation of Rotational Bracket on a Link $L_5$}
     \label{linkcalc5}
\end{center}
\end{figure}

In Figures~\ref{linkcalc1},~\ref{linkcalc2},~\ref{linkcalc3},~\ref{linkcalc4} and ~\ref{linkcalc5} we calculate the bracket for links respectively labeled $L_{1}, L_{2}, L_{3}, L_{4}, L_{5}.$\\

$$ <L_{1}> = <L_{2}> =  -d^{2} + \DegZero^{2} + \DegZero \DegTwo.$$
$$<L_{3}> = d \DegZero.$$
$$<L_{4}> = d^{2}.$$

While $L_{1}$ and $L_{2}$ have the same bracket invariant, they are distinct links. The calculations show that they are both non-trivial.
Neither $L_{3}$ nor $L_{4}$ are distinguished from their corresponding split links but neither of them is a split rotational link.
(These results are proved in Section 3.3 by using the techniques of the parity bracket described in that section.)
The result for $<L_{5}>$ is indicated in Figure~\ref{linkcalc5}. This calculation verifies the non-triviality of $L_{5}.$ More information about these links may come from other quantum link invariants.
We will investigate these examples more deeply in a sequel to this paper.\\

\subsection{Parity Bracket Polynomial for Rotational Virtuals}

\begin{figure}
     \begin{center}
     \begin{tabular}{c}
     \includegraphics[width=6cm]{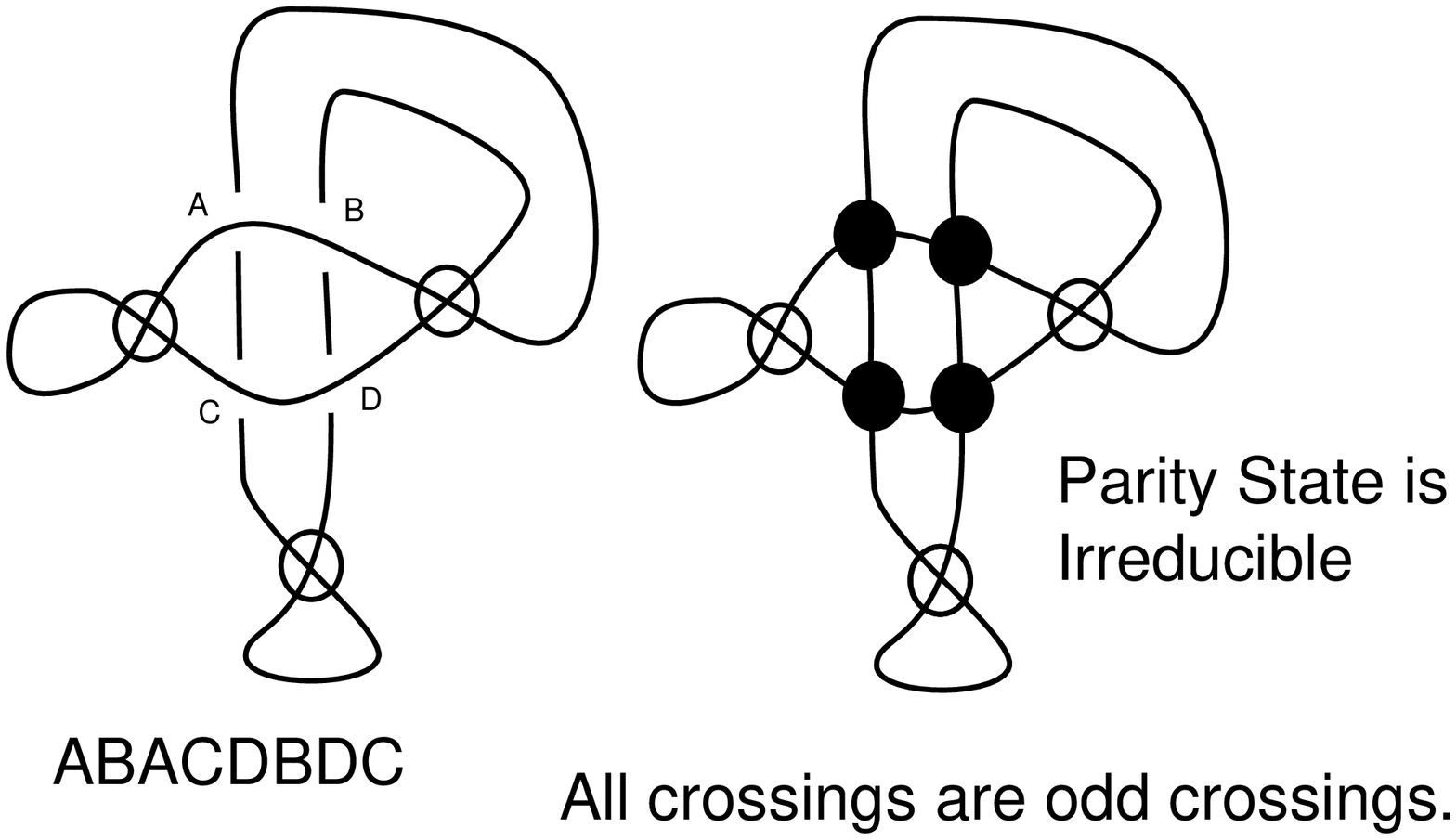}
     \end{tabular}
     \caption{\bf Single State Detection of a Rotational Virtual Knot}
     \label{parity}
\end{center}
\end{figure}

We now generalize the Parity Bracket Polynomial of Vassily Manturov \cite{MP,Ma,MP1,PCFree,CobFree} to rotational virtuals. 
This original parity bracket is a generalization of the bracket polynomial to virtual knots and links that uses the parity of the crossings. In the discussion below, we give the definition of the parity bracket for rotational virtuals. The definition is obtained by using the regular homotopy detour move as in the case of the rotational bracket polynomial. Thus there are end--states of the rotational bracket that have virtual crossings and these can now occur in graphs with $4$-valent nodes immersed in the plane or on the $2$--sphere. \\

Recall that a non-virtual crossing  $X$ in a diagram is said to be {\it even} if, in the Gauss word for the diagram, there are an even number of crossing labels (non-virtual) appearing between the two occurrences of $X$ in the code. A crossing with an odd number of such occurrences is said to be {\it odd}. See Figure~\ref{parity} for an example. In that diagram the code word is 
$ABACDBDC$ and we see that every crossing is odd.  From now on, when we refer to crossings it will be assumed that they are not virtual crossings unless so specified. For another example, see
Figure~\ref{Figure 12} where there are four odd crossings and one even crossing.\\

\begin{figure}
     \begin{center}
     \begin{tabular}{c}
     \includegraphics[width=8cm]{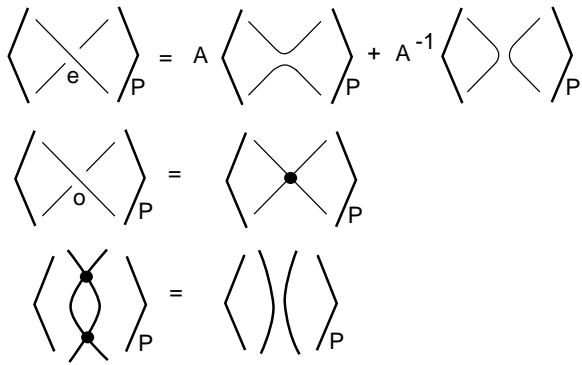}
     \end{tabular}
     \caption{\bf Parity Bracket Expansion}
     \label{Figure 10}
\end{center}
\end{figure}

\begin{figure}
     \begin{center}
     \begin{tabular}{c}
     \includegraphics[width=6cm]{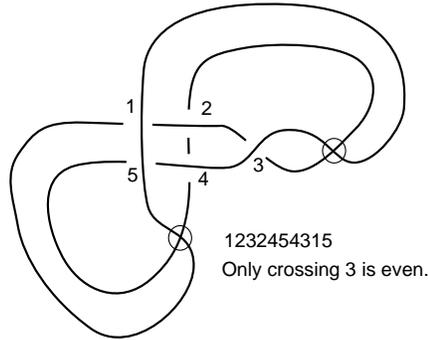}
     \end{tabular}
     \caption{\bf A Knot KS With Unit Jones Polynomial}
     \label{Figure 12}
\end{center}
\end{figure}

We define a {\it rotational virtual graph} to be a locally $4$--regular graph that is immersed in the plane or $2$--sphere with multiple points for the immersion corresponding to virtual crossings in the graph. Two such graphs are {\it rotationally equivalent} if one can be obtained from the other by planar (or $2$-sphere) isotopy coupled with regular homotopy detour moves for consecutive sequences of virtual crossings. A rotational virtual graph is said to be {\it reduced} if there are no simplifying graphical Reidemeister two-moves available on it. The type two graphical move is illustrated in the third diagram of 
Figure~\ref{Figure 10}. Note that the presence of such a move is a combinatorial condition that can be deduced from the immersion of the graph in the plane. In practice, it may be useful to apply detour moves to exhibit the possibility of a reduction. Since we are dealing with rotational virtuals here, a virtual curl on an edge in the graph can obstruct a reduction. For example, the graph shown in 
Figure~\ref{parity} is irreducible.\\
 
We define a {\em parity state} of a rotational virtual diagram $K$ to be a labeled rotational virtual
graph obtained from $K$ as follows: For each odd crossing in $K$ replace the crossing by a graphical node. For each even crossing in $K$ replace the crossing by one of its two possible smoothings, and label the smoothing site by $A$ or $A^{-1}$ in the usual way. Then we define the parity bracket by the 
state expansion formula
$$\langle K \rangle _{P} = \sum_{S}A^{n(S)}[S]$$
where $n(S)$ denotes the number of $A$-smoothings minus the number of $A^{-1}$ smoothings and 
$[S]$ denotes a combinatorial evaluation of the state defined as follows: In this expansion reduce each state
by graphical Reidemeister two-moves on nodes as shown in Figure~\ref{Figure 10} and by regular homotopy. The graphs are taken up to virtual equivalence (planar isotopy plus detour moves on the virtual crossings). Then regard the reduced state as a disjoint union of standard state loops (without nodes but these may have virtual crossings) and graphs that irreducibly contain nodes and possibly virtual crossings.
With this we write $$[S] = (- A^{2} - A^{-2})^{l(S)} [G(S)]$$ where $l(S)$ is the number of standard circular loops in the reduction of the state $S$ and $[G(S)]$ is the disjoint union of reduced graphs that contain nodes and/or virtual crossings.\\

In this way, we obtain a sum of Laurent polynomials in $A$ multiplying reduced graphs as the {\it rotational parity bracket.} It is not hard to see that this bracket is invariant under isotopy and rotational detour moves and that it behaves just like the usual bracket under the first Reidemeister move. However, the use of
parity to make this bracket expand to graphical states gives it considerable extra power in some situations. For example, consider the diagram in Figure~\ref{parity}. We see that all the classical crossings in this knot are odd. Thus the parity bracket is just the graph obtained by putting nodes at each of these crossings. The resulting graph does not reduce under the graphical two moves, and so we conclude that this rotational knot  is non-trivial and non-classical. Since we can apply the parity
bracket to a flat knot by taking $A = -1$, we see that this method shows that the  flat rotational version of this knot is non-trivial. \\

\begin{figure}
     \begin{center}
     \begin{tabular}{c}
     \includegraphics[width=6cm]{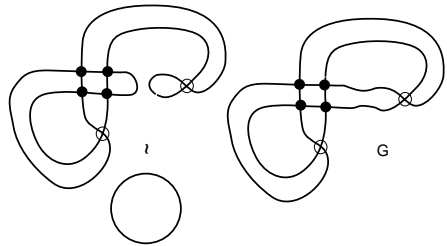}
     \end{tabular}
     \caption{\bf Parity Bracket States for the Knot KS}
     \label{Figure 13}
\end{center}
\end{figure}

In Figure~\ref{Figure 12} we show an example knot $KS$ and its parity state expansion in in Figure~\ref{Figure 13}.  The first state reduces to a circle. 
The second graphical state is irreducible, proving that $KS$ is a non-trivial rotational virtual knot.\\

We end this subsection with a theorem that shows, in principle, how powerful is this use of parity.\\

\noindent {\bf Irreducibility Theorem.} Let $D$ be a rotational virtual knot or link diagram with all odd crossings. Then there is a diagram $D',$ obtained from $D$ by adding local flat virtual curls at some arcs of the 
diagram $D,$ such that the rotational virtual graph $G(D')$ obtained from $D'$ by replacing all crossings in $D'$ by graphical nodes is irreducible in the sense above. Thus the rotational virtual link $D'$ is non-trivial. Furthermore, the rotational genus $Rg(D')$ is equal to $g(D'),$ the genus of the diagram $D',$ as defined in Section 3.1.\\

\begin{figure}[htb]
     \begin{center}
     \begin{tabular}{c}
     \includegraphics[width=5cm]{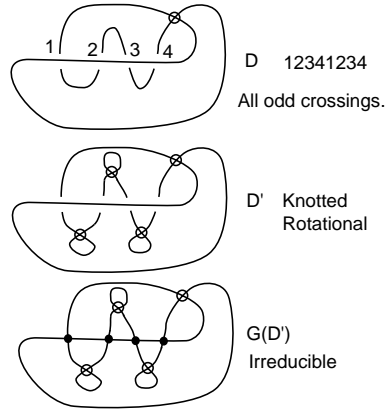}
     \end{tabular}
     \caption{\bf All Odd Placement}
     \label{placement}
\end{center}
\end{figure}

\begin{figure}[htb]
     \begin{center}
     \begin{tabular}{c}
     \includegraphics[width=5cm]{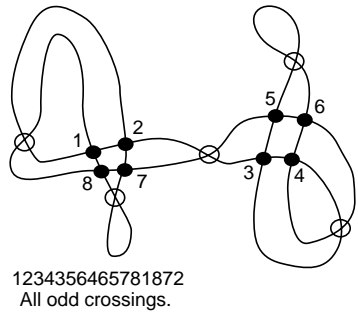}
     \end{tabular}
     \caption{\bf Genus Two}
     \label{genustwo}
\end{center}
\end{figure}

\begin{figure}
     \begin{center}
     \begin{tabular}{c}
     \includegraphics[width=6cm]{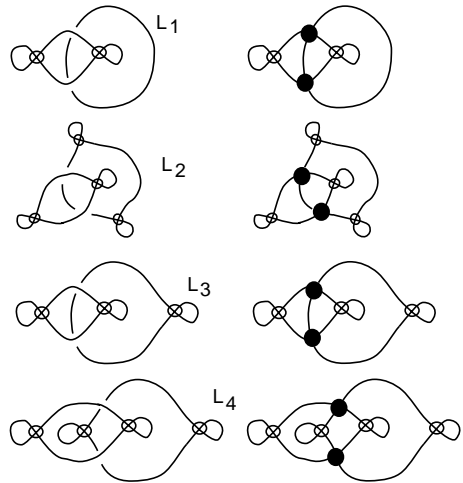}
     \end{tabular}
     \caption{\bf Parity Link Bracket proves that $L_{1},L_{2},L_{3},L_{4}$ are all non-trivial and pairwise distinct.}
     \label{graphlink}
\end{center}
\end{figure}

\noindent{\bf Proof.} Note that by placing virtual flat curls at appropriate places on the diagram $D$ to produce a diagram $D'$, we can ensure that the graph $G(D')$ can not be reduced by any type two  moves. This means that the genus for the diagram $D'$ is minimal among all diagrams rotationally equivalent to $D'.$  This completes the proof. // \\

\noindent {\bf Remark.} The example in Figure~\ref{parity} can be regarded as an example of this theorem applied to a diagram $D$ that does not have the flat virtual curl between crossings $A$ and $C.$ In Figure~\ref{placement} we give an explicit example of placing curls to create an irreducible diagram. In Figure~\ref{genustwo} we give an irreducible graph such that any virtual diagram that overlies it (replacing graphical nodes by classical crossings) will have this graph as its parity bracket. In this case the reader can verify that the genus of the diagram is equal to two.\\

\noindent{\bf Remark.} The parity invariant generalizes for links by replacing crossings between components by graphical nodes (and continuing to reduce the graphs by graphical Reidemeister Two moves).
Thus we find as shown in Figure~\ref{graphlink} That all the examples $L_{1}, L_{2}, L_{3}$ and $L_{4}$ of the previous section are non-split links and are distinct from one another in the rotational category.
See our paper with Aaron Kaestner \cite{KaestnerKauffman} for more about parity and links.\\

Parity is clearly an important theme in virtual knot theory and will figure in many future investigations of this subject. The type of construction that we have indicated for the bracket polynomial in this section can be varied and applied to other invariants. Furthermore the notion of describing a parity for crossings
in a diagram is also susceptible to generalization. For more on this theme the reader should consult
\cite{MP1} and \cite{SL} for our original use of parity for another variant of the bracket polynomial. The irreducibility theorem illustrates the striking power of parity in that we can use it to find infinitely many examples of rotational virtual diagrams that are essentially irreducible in the sense discussed above.\\

In the next section we will see the general formulation of quantum link invariants and how the bracket invariant of this section specializes to certain matrix models using cup, cap and Yang-Baxter operators.\\

\section{Quantum Link Invariants}

There are virtual link invariants corresponding to every quantum link invariant
of classical links. However this must be said with a caveat: We do not assume
invariance under the first classical Reidemeister move (hence these are
invariants of regular isotopy) and we do not assume invariance under the flat
version of the first Reidemeister move in the (B) list of virtual moves. 
Otherwise the usual tensor or state sum formulas for quantum link invariants
extend to this generalized notion of regular isotopy invariants of virtual knots
and links.\\

\begin{figure}
     \begin{center}
     \begin{tabular}{c}
     \includegraphics[width=10cm]{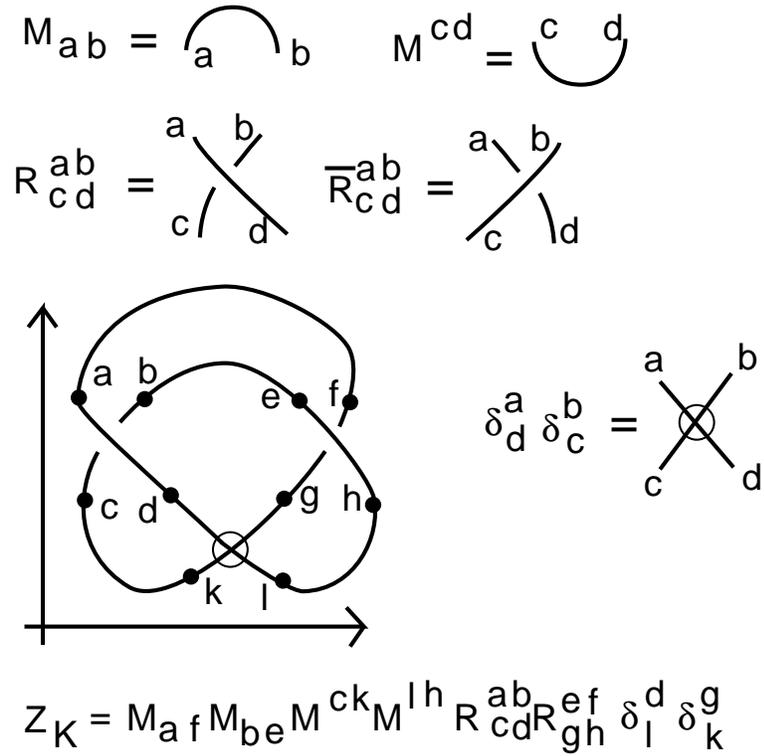}
     \end{tabular}
     \caption{\bf Quantum Link Invariants}
     \label{fig20}
\end{center}
\end{figure}

In order to carry out this program, we quickly recall how to construct quantum
link invariants in the unoriented case. See \cite{KP} for more details.  The abstract definition of a quantum link invariant is usually formulated via a functor from the Tangle Category 
to the a category of modules over a commutative ring with unit. Here we use the matrix forms that are the end result of such a formulation. Thus we build a ``statistical mechanics" model by associating matrices to the cup, caps and crossings of a Morse diagram, a diagram equipped with a height function. To this purpose the link
diagram is arranged with respect to a given ``vertical" direction in the plane so
that perpendicular lines to this direction intersect the diagram transversely or
tangentially at maxima and minima.. In this way the diagram can be seen as
constructed from a pattern of interconnected maxima, minima and crossings---as
illustrated in Figure~\ref{fig20}.\\

\begin{figure}
     \begin{center}
     \begin{tabular}{c}
     \includegraphics[width=10cm]{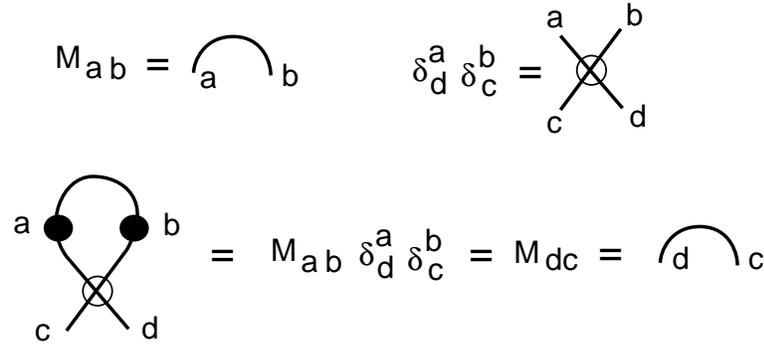}
     \end{tabular}
     \caption{\bf Quantum Virtual Curl}
     \label{vcurl}
\end{center}
\end{figure}

As illustrated in Figure~\ref{fig20}, we associate symbols $M_{ab}$ and $M^{ab}$ to 
maxima and minima respectively, and symbols $R^{ab}_{cd}$ and $\overline{R}^{ab}_{cd}$
to the two types of crossings. The indices on these symbols indicate how they are
interconnected. Each maximum or minimum has two lines available for connection
corresponding to the indices $a$ and $b.$ Each $R$ , $\overline{R}$  has four
lines available for connection.  Thus the symbol sequence

$$T(K) = M_{af}M_{be}M^{ck}M^{lh}R^{ab}_{cd}R^{ef}_{gh}\delta^{d}_{l}\delta^{g}_{k}$$

\noindent represents the virtual trefoil knot as shown in Figure~\ref{fig20}.  Since repeated
indices show the places of connection, there is no necessary order for this
sequence of symbols. I call $T(K)$ an {\em abstract tensor} expression for the
virtual trefoil knot $K$. \\

By taking matrices (with entries in a commutative ring) for the $M$'s and the
$R$'s it is possible to re-interpret the abstract tensor expression as a
summation of products of matrix entries over all possible choices of indices in
the expression.  Appropriate choices of matrices give rise to link invariants. 
If $K$ is a knot or link and $T(K)$ its associated tensor expression, let $Z(K)$
denote the evaluation of the tensor expression that corresponds to the above
choice of matrices. We will assume that the matrices have been chosen so that
$Z(K)$ is an invariant of regular isotopy. \\

Here is the list of the algebraic versions of the topological moves. Move $0$ is
the cancellation of maxima and minima. Move $II$  corresponds to the second
Reidemeister move. Move $III$ is the Yang-Baxter equation.  Moves $IV$ and $IV'$ express
the relationship of switching a line across a maximum. (There is a corresponding
version where the line is switched across a minimum, but it is a consequence of the other axioms.) Implicit in this list of moves are the 
formal identities for manipulating Kronecker deltas such as $\delta^{a}_{b} M^{bc} = M^{ac}$ and 
$\delta^{a}_{b}\delta^{b}_{c} = \delta^{a}_{c}.$ Note that in working with abstract tensors, one must obey the usual substitution rules for indices that one
uses in standard tensor calculus. That is, if an index is replaced at any place in an expression, it must be replaced at all occurrences of that index. The replacement should be 
a symbol different from other indices used in the given expression.\\

$$0. \,\,  \hspace{.1in} M^{ai}M_{ib}  =\delta^{a}_{b}$$
$$II.\hspace{.1in} R^{ab}_{ij} \overline{R^{ij}_{cd}} = \delta^{a}_{c}
\delta^{b}_{d}$$
$$III.\hspace{.1in} R^{ab}_{ij}R^{jc}_{kf}R^{ik}_{de} =
R^{bc}_{ij}R^{ai}_{dk}R^{kj}_{ef}$$
$$IV. \hspace{.1in}R^{ai}_{bc}M_{id} = M_{bi}  \overline{R^{ia}_{cd}}$$
$$IV'. \hspace{.1in}\overline{R^{ai}_{bc}}M_{id} = M_{bi}  R^{ia}_{cd}$$

\begin{figure}
     \begin{center}
     \begin{tabular}{c}
     \includegraphics[width=5cm]{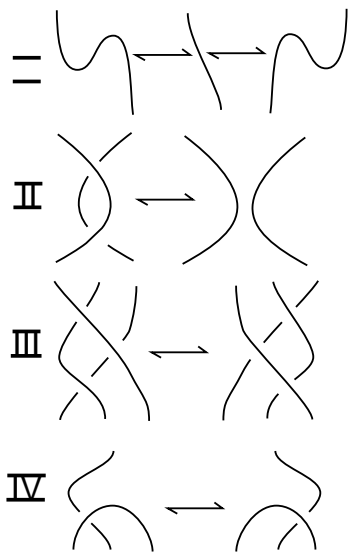}
     \end{tabular}
     \caption{\bf Regular Isotopy with respect to a Vertical Direction}
     \label{regvert}
\end{center}
\end{figure}

\noindent These equations encapsulate the abstract structure of unoriented quantum link invariants, and correspond to invariance under the moves shown in 
Figure~\ref{regvert}. This reader may enjoy taking each part of this figure (we have drawn representative cases) and comparing it with the corresponding abstract tensor equation above.\\

The generalisation of the quantum link invariant $Z(K)$ to virtual knots and
links is quite straightforward. {\it We simply ignore the virtual crossings in the
diagram.} Another way to put this is that we take each virtual crossing to be
represented by crossed Kronecker deltas as in Figure~\ref{fig20}.  The virtual crossing is
represented by the tensor  $$V^{ab}_{cd} = \delta^{a}_{d} \delta^{b}_{c}.$$ Here
$\delta^{a}_{b}$ is the Kronecker delta. In a matrix representation $\delta^{a}_{b}$ is equal to $1$ if $a=b$ and is equal
to $0$ otherwise. (The Kronecker delta is well defined as an abstract
tensor by its formal properties.) \\

In extending $Z(K)$ to virtual knots and links  by this method we cannot hope to
obtain invariance under the type I virtual move. In fact, as Figure~\ref{vcurl} shows, the
presence of a virtual curl is indexed by the transpose $M_{dc}$ of the tensor
$M_{cd}$. Thus we have defined {\em virtual rotational isotopy} to be invariance under all
the extended Reidemeister moves for virtuals except type (A)I and (B)I.  It is
easy to see that $Z(K)$ extends in this way when $Z(K)$ is an invariant of
rotational isotopy for classical links. \\

In particular the bracket polynomial for classical knots and links is obtained by letting
the indices run over the set $\{1,2\}$ with $M^{ab} = M_{ab}$ for all $a$ and $b$
and $M_{11}=M_{22} = 0$ while $M_{12} = iA$ and $M_{21} = -iA^{-1}$ where $i^{2}
= -1.$  Thus 
$$M =  
\left[ \begin{array}{cc} 
0& iA  \\ 
-iA^{-1} & 0 
\end{array} 
\right]. 
$$

\noindent The $R$ matrices are defined by the equations

$$R^{ab}_{cd} = A M^{ab}M_{cd} + A^{-1} \delta^{a}_{c} \delta^{b}_{d}$$

$$\overline{R}^{ab}_{cd} = A^{-1} M^{ab}M_{cd} + A \delta^{a}_{c}
\delta^{b}_{d}$$

\noindent These equations for the $R$ matrices are the algebraic translation of the
smoothing identities for the bracket polynomial.    Then we have\\

\noindent {\bf Theorem}  With $Z(K)$ defined as above and $K$ a classical knot
or link, then $$Z(K) = d<K>$$ where $d = -A^{2} - A^{-2}$  and $<K>$ denotes the topological bracket polynomial.\\

\noindent {\bf Proof.} See \cite{KP}.// \\

For the extension of $Z(K)$  to virtuals there is a state summation similar to
that of the bracket polynomial. For this, let $C$ be a diagram in the plane that
has only virtual crossings. View this diagram as an immersion of a circle in the
plane. Let $rot(C)$ denote the absolute value of the Whitney degree of $C$ as a
immersion in the plane. (Since $C$ is unoriented only the absolute value of the
Whitney degree is well-defined.). The Whitney degree of an oriented plane
immersion is the total algebraic number of $2 \pi$ turns of the unit tangent
vector to the curve as the curve is traversed once. Let $d(C)$ be defined by the
equation $$d(C) = (-1)^{rot(C)} (A^{2 rot(C)} + A^{-2 rot(C)}).$$ Let $S$ be a
state of a virtual diagram $K$ obtained by smoothing each classical crossing in
$K.$  Let $C \in S$ mean that $C$ is one of the curves in $S.$  Let $<K|S>$
denote the usual product of vertex weights ($A$ or $A^{-1}$) in the bracket state
sum.  Then \vspace{3mm}

\noindent {\bf Proposition}  The invariant of virtual regular isotopy $Z(K)$
is described by the following state summation.

$$Z(K) = \sum_{S} <K|S> \prod_{C \in S} d(C)$$

\noindent where the terms in this formula are as defined above. Note that $Z(K)$
reduces to $d<K>$ when $K$ is a classical diagram. \vspace{3mm}

\noindent {\bf Proof.} The proof is a calculation based on the tensor model
explained in this section. The details of this calculation are omitted.//\\

\noindent {\bf Remark.}  The state sum in this Proposition  generalizes to an
invariant of virtual regular isotopy with an infinite number of polynomial
variables, one for each regular homotopy class of unoriented curve immersed in
the $2$--sphere. To make this generalization, let $d_{n}$ for $n=0,1,2,3,...$ denote a
denumerable set of commuting independent variables.  If $C$ is an immersed curve
in the plane, define $Var(C) = d_{n}$ where $n = rot(C)$, the absolute value of
the Whitney degree of $C.$  We take $d_{1} = -A^{2} - A^{-2}$ as before, but the
other variables are independent of each other and  of $A$. \\

\noindent Now define the generalisation of $Z(K)$, denoted $\overline{Z}(K),$ by
the formula

$$\overline{Z}(K) = \sum_{S} <K|S> \prod_{C \in S} Var(C).$$

In this definition we have replaced the evaluation $d(C)$ by the corresponding
variable $Var(C).$  The reader will note that this generalization is exactly the rotational version of the bracket that we have studied in the previous section of this paper.\\

\subsection{The Binary Bracket Polynomial as a Quantum Link Invariant}
The reader will easily verify that the following matrix is a solution to the Yang-Baxter Equation and that for unoriented models, as discussed above, the cup and cap operators can be taken to be the identity matrix.

$$R = \left( \begin{array}{cccc}
0  & 0  & 0  & A \\
0  & A^{-1} & 0 & 0 \\
0 & 0  & A^{-1} & 0 \\
A & 0 & 0 & 0 \\
\end{array} \right).$$

This $4 \times 4$ matrix is viewed as acting upon a tensor product of a two-dimensional space with itself whose basis indices are
$0$ and $1.$ Note that if $A$ is a unit complex number, then $R$ is a unitary matrix. This makes this matrix of interest in the context of quantum computing as well as topology. See \cite{Dye,BraidGates}.\\

The corresponding quantum link invariant (for arbitrary $A$) satisfies the special bracket expansion that we have elsewhere \cite{SL} called the {\it binary bracket polynomial}. Here the resulting invariant
of virtual knots and links is also invariant under the first flat virtual move and so it does not give new information about rotational virtuals. Nevertheless, the invariant is of interest.\\

We now describe the binary bracket as a state summation. In this respect, it has almost exactly the same formalism as
the standard bracket polynomial, except that the value of an unlabeled loop is equal to $2,$ and the loops in each state are
colored with the  colors from the set $\{0,1 \}$ in such a way that the colors appearing at a smoothing are always different.
This restricts the possible states to a very small number and causes the invariant to behave differently on virtual links than
it does on classical links.
\bigbreak

Let $K$ be any unoriented (virtual) link diagram. Define an {\em unlabeled state}, $S$, of $K$  to
be a choice of smoothing for each  classical crossing of $K.$ There are two choices for smoothing a given  crossing, and
thus there are $2^{N}$ unlabeled states of a diagram with $N$ classical crossings.  A {\em labeled state} is a state $S$ such
that the labels $0$ (zero) or $1$ (one) have been assigned to each component loop in the state. 
\bigbreak

\begin{figure}
     \begin{center}
     \begin{tabular}{c}
     \includegraphics[width=8cm]{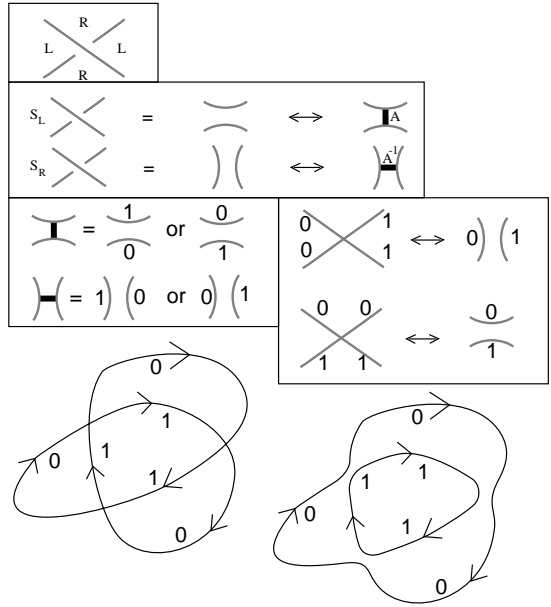}
     \end{tabular}
     \caption{\bf Bracket Smoothings}
     \label{bsmooth}
\end{center}
\end{figure}

In  a state we designate each smoothing with $A$ or $A^{-1}$ according to the left-right convention 
shown in Figure~\ref{bsmooth}. This designation is called a {\em vertex weight} of the state. {\em We require of a labeled state
that the two labels that occur at a smoothing of a crossing are distinct.} This is indicated by a bold line between the 
arcs of the smoothing as illustrated in Figure~\ref{bsmooth}.  
Labeled states satisfying this condition at the site of every smoothing will be called {\em properly labeled states.}
If $S$ is a properly labeled state, we let $\{ K|S \}$ denote the product of its vertex weights, and we define the 
two-color bracket polynomial by the equation:   
$$\{ K \} \, = \sum_{S} \{ K|S \}.$$
where $S$ runs through the set of properly labeled states of $K.$\\

It follows from this definition that $\{ K \}$ satisfies the equations
 $$\includegraphics[width=5cm]{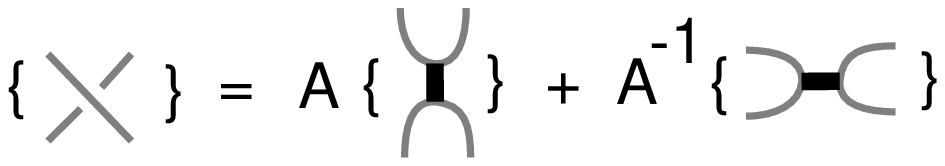}$$ 
$$\{ K \amalg  O \} \, = 2 \{ K \},$$
$$\{ O \} \, =2.$$
  
\noindent The first equation expresses the fact that the entire set of states of a given diagram is
the union, with respect to a given crossing, of those states with an $A$-type smoothing and those
 with an $A^{-1}$-type smoothing at that crossing. In the first equation, we indicate that the colors at the smoothing are
different by the dark band placed between the arcs of the smoothing. The second and the third equations are clear from the formula
defining the state summation. See Figure~\ref{bsmooth} for a summary of the conventions for coloring states of the binary bracket polynomial.
\bigbreak
 
The {\em  binary bracket polynomial} , $\{ K \} \, = \, \{ K \}(A)$,  assigns to each unoriented (virtual) link diagram $K$ a 
Laurent polynomial in the variable $A.$
\bigbreak

In computing the  binary bracket, one finds the following behaviour under Reidemeister move I: 
$$\{\Rcurl \} = A \{ \Arc \},$$ 
$$\{ \Lcurl \} = A^{-1} \{ \Arc \}. $$

\noindent {\bf Theorem.} The  binary bracket is invariant under regular isotopy for virtual links, and it can be  normalized to
an invariant of ambient isotopy by the definition  
$$Inv_{K}(A) = A^{-w(K)} \{K \}(A),$$ where we choose an orientation for $K$, and where $w(K)$ is 
the sum of the crossing signs  of the oriented link $K$. $w(K)$ is called the {\em writhe} of $K$. 
\bigbreak

\noindent{\bf Discussion.} We have given the basic construction for the state summation for the binary bracket polynomial for virtual knots.
While the quantum link invariant version of the binary bracket is invariant under the addition or removal of flat virtual curls, the {\it state summation} for the binary bracket
can be generalized to an invariant of rotational virtuals just as we have done in the previous section for the bracket polynomial. We leave out the details and examples for this construction, but encourage the reader to explore this avenue. Many interesting phenomena arise and they will be explored in a separate paper.\\

\section{Oriented Quantum Link Invariants} 

Slight but significant modifications are needed
to write the oriented version of the models we have discussed in the previous
section. See \cite{KP}, \cite{Turaev}, \cite{LOMI}, \cite{Hennings}. In this
section we sketch the construction of oriented topological amplitudes.\\

We include this section of the paper as a basis for further explorations of rotational virtual link invariants. Essentially all quantum link invariants fit into the framework described in this section.
We use this description to show the well-known quantum link invariant models for specializations of the Homflypt polynomial and their rotational generalizations. We will continue with more
examples in a subsequent paper. Here the reader will find a concise description of the structure of these invariants.\\

The generalization  to oriented link diagrams naturally involves the introduction
of right and left oriented caps and cups. These are drawn as shown in Figure~\ref{conversion}
below.\\

\begin{figure}
     \begin{center}
     \begin{tabular}{c}
     \includegraphics[width=5cm]{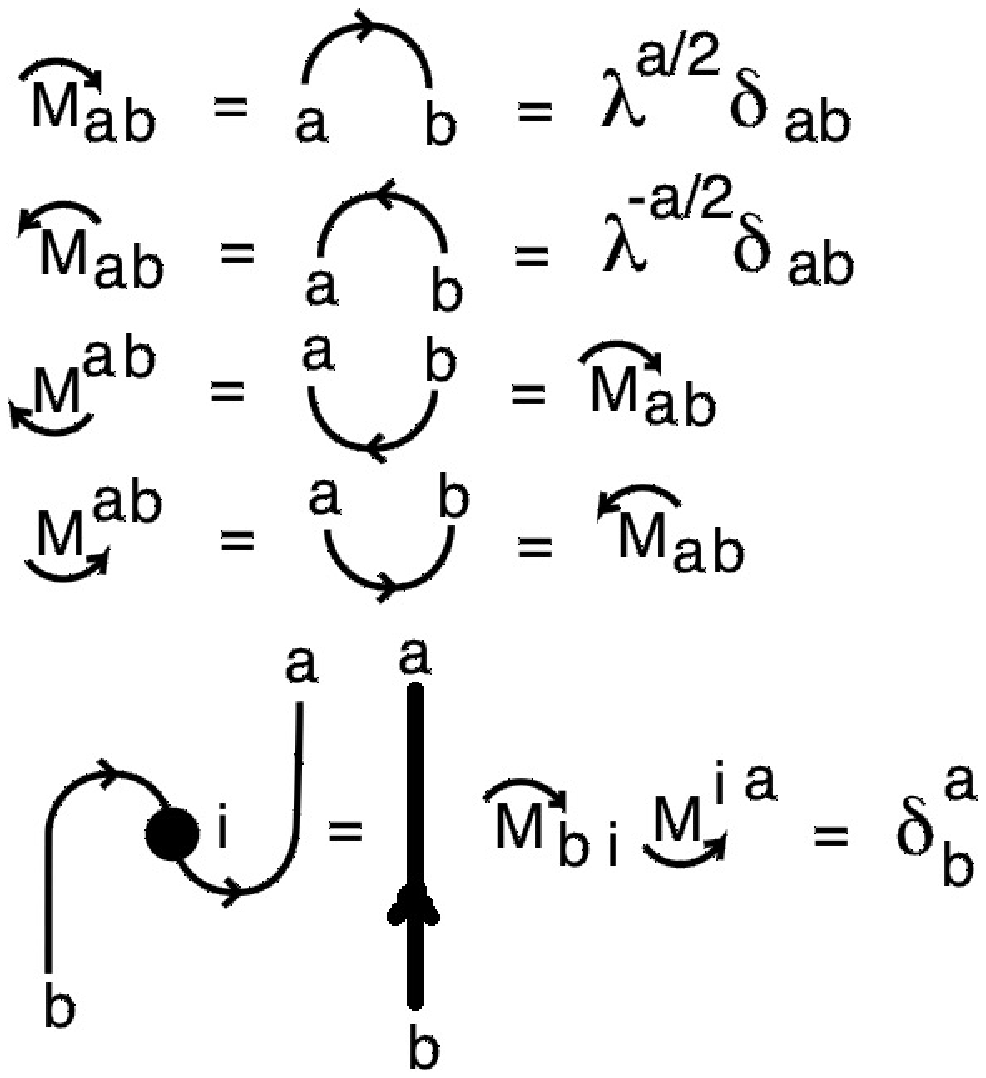}
     \end{tabular}
     \caption{\bf Right and Left Cups and Caps}
     \label{rlcupcap}
\end{center}
\end{figure}

A right cup cancels  with a right cap to produce an upward pointing identity
line.  A left cup cancels with a left cap to produce a downward pointing identity
line. \vspace{3mm}

Just as we considered the simplifications that occur in the unoriented model by
taking the cup and cap matrices to be identical, lets assume here that right caps
are identical with left cups and that consequently left caps are identical with
right cups.  In fact, let us assume that the right cap and left cup are given by
the matrix $$M_{ab} = \lambda^{a/2} \delta_{ab} = M^{ab}$$ where $\lambda$ is a constant
to be determined by the situation, and $\delta_{ab}$ denotes the Kronecker delta.
 Then the left cap and right cup are given by the inverse of $M$:
$$M_{ab}^{-1} =
\lambda^{-a/2} \delta_{ab} = (M^{ab})^{-1}.$$

We assume that along with $M$ we are given a solution $R$ to the Yang-Baxter
equation, and that in an oriented diagram the specific choice of $R^{ab}_{cd}$ is
governed by the local orientation of the crossing in the diagram. Thus $a$ and
$b$ are the labels on the lines going into the crossing and $c$ and $d$ are the
labels on the lines emanating from the crossing.
\newpage

\begin{figure}
     \begin{center}
     \begin{tabular}{c}
     \includegraphics[width=5cm]{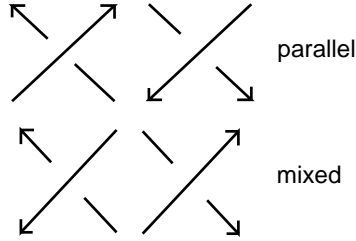}
     \end{tabular}
     \caption{\bf Oriented Crossings}
     \label{orcross}
\end{center}
\end{figure}

Note that with respect to the vertical direction for the amplitude, the crossings
can assume the aspects: both lines pointing upward, both lines pointing downward,
one line up and one line down (two cases). See Figure~\ref{orcross}.\\

Call the cases of one line up and one line down the {\em mixed} cases and the
upward and downward cases the {\em parallel} cases. A given mixed crossing can be
converted ,in two ways, into a combination of a parallel crossing of the same
sign plus a cup and a cap. See Figure~\ref{conversion}.\\

\begin{figure}
     \begin{center}
     \begin{tabular}{c}
     \includegraphics[width=5cm]{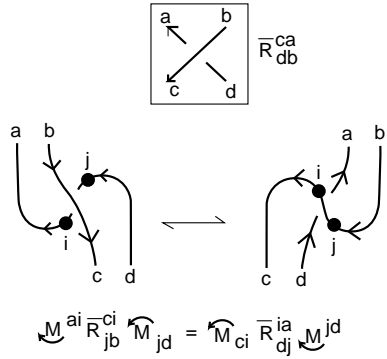}
     \end{tabular}
     \caption{\bf Conversion}
     \label{conversion}
\end{center}
\end{figure}

This leads to an equation that must  be satisfied by the $R$ matrix in relation
to powers of $\lambda$ (again we use the Einstein summation convention):

$$\lambda^{a/2} \delta^{ai} R^{ci}_{jb} \lambda^{-d/2} \delta_{jd} =
\lambda^{-c/2} \delta_{ic} R^{ia}_{dj} \lambda^{b/2} \delta^{jb}.$$

This simplifies to the equation

$$\lambda^{a/2}   R^{ca}_{db} \lambda^{-d/2}  = \lambda^{-c/2} R^{ca}_{db}
\lambda^{b/2} ,$$

\begin{figure}
     \begin{center}
     \begin{tabular}{c}
     \includegraphics[width=5cm]{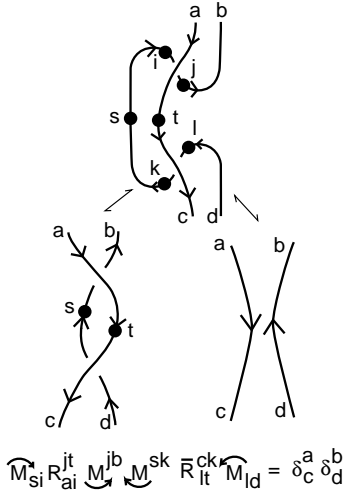}
     \end{tabular}
     \caption{\bf Antiparallel Second Move}
     \label{anti}
\end{center}
\end{figure}

\noindent from which we see that $R^{ca}_{db}$ is necessarily equal to zero unless $b+d =
a+c$.  We say that the $R$ matrix is {\em spin preserving} when it satisfies this
condition.  Assuming that the $R$ matrix is spin preserving,  the model will be
invariant under all orientations of the second and third Reidemeister moves just
so long as it is invariant under the anti-parallel version of the second
Reidemeister move as shown in Figure~\ref{anti}.\\

This antiparallel version of the second Reidemeister move places the following
demand on the relation between $\lambda$ and $R$:

$$\sum_{st} \lambda^{(s-b)/2} \lambda^{(t-c)/2}R^{bt}_{as}
\overline{R^{cs}_{dt}}= \delta^{a}_{c} \delta^{b}_{d}.$$

Call this the {\em $R-\lambda$ equation.} The reader familiar with \cite{JO}
or with the piecewise linear version as described in \cite{KP} will
recognise this equation as the requirement for regular isotopy invariance in
these models. \vspace {3mm}

\noindent{\bf Remarks and An Example}. Note that oriented invariants as described here will not, in general, be invariant under the flat virtual curl move since we have
that the right cap and left cup are given by
the matrix $$M_{ab} = \lambda^{a/2} \delta_{ab},$$  and the left cap and right cup are given by the inverse of $M$:
$$M_{ab}^{-1} = \lambda^{-a/2} \delta_{ab}.$$ See Figure~\ref{rlcupcap}.
This means that, using our usual extension to rotational invariants by using a permutation operator for virtual crossings will yield non-trivial rotational 
invariants, just as in the case of the matrix model for the bracket polynomial.\\

\begin{figure}
     \begin{center}
     \begin{tabular}{c}
 {\tt    \setlength{\unitlength}{0.92pt}
\begin{picture}(428,232)
\thinlines    \put(338,47){$D'$}
              \put(262,48){$C'$}
              \put(191,49){$B'$}
              \put(319,174){$D$}
              \put(254,175){$C$}
              \put(181,175){$B$}
              \put(347,191){\framebox(38,30){$\ne$}}
              \put(358,69){\framebox(38,30){$\ne$}}
              \put(16,137){\vector(3,4){63}}
              \put(86,142){\vector(-1,1){33}}
              \put(43,185){\vector(-1,1){33}}
              \put(349,11){\vector(3,4){63}}
              \put(417,15){\vector(-1,1){76}}
              \put(123,140){\vector(0,1){77}}
              \put(155,140){\vector(0,1){78}}
              \put(225,139){\vector(0,1){78}}
              \put(194,140){\vector(0,1){77}}
              \put(293,139){\vector(0,1){78}}
              \put(265,139){\vector(0,1){77}}
              \put(90,22){\vector(-1,1){76}}
              \put(18,22){\vector(1,1){30}}
              \put(59,63){\vector(1,1){34}}
              \put(281,14){\vector(0,1){77}}
              \put(309,14){\vector(0,1){78}}
              \put(210,14){\vector(0,1){77}}
              \put(238,13){\vector(0,1){78}}
              \put(162,13){\vector(0,1){78}}
              \put(133,13){\vector(0,1){77}}
              \put(404,140){\vector(-1,1){76}}
              \put(336,136){\vector(3,4){63}}
              \put(199,170){\framebox(21,20){$<$}}
              \put(270,168){\framebox(19,21){$>$}}
              \put(138,41){\framebox(20,19){$=$}}
              \put(130,167){\framebox(20,19){$=$}}
              \put(214,43){\framebox(21,20){$<$}}
              \put(285,42){\framebox(19,21){$>$}}
              \put(80,165){\makebox(38,23){$= A$}}
              \put(160,166){\makebox(25,23){$+$}}
              \put(231,165){\makebox(21,23){$+$}}
              \put(297,164){\makebox(27,24){$+$}}
              \put(103,42){\makebox(27,21){$= A'$}}
              \put(166,39){\makebox(26,25){$+$}}
              \put(243,39){\makebox(23,22){$+$}}
              \put(313,37){\makebox(30,26){$+$}}
\end{picture}}
     \end{tabular}
     \caption{\bf Local States}
     \label{homflypt}
\end{center}
\end{figure}
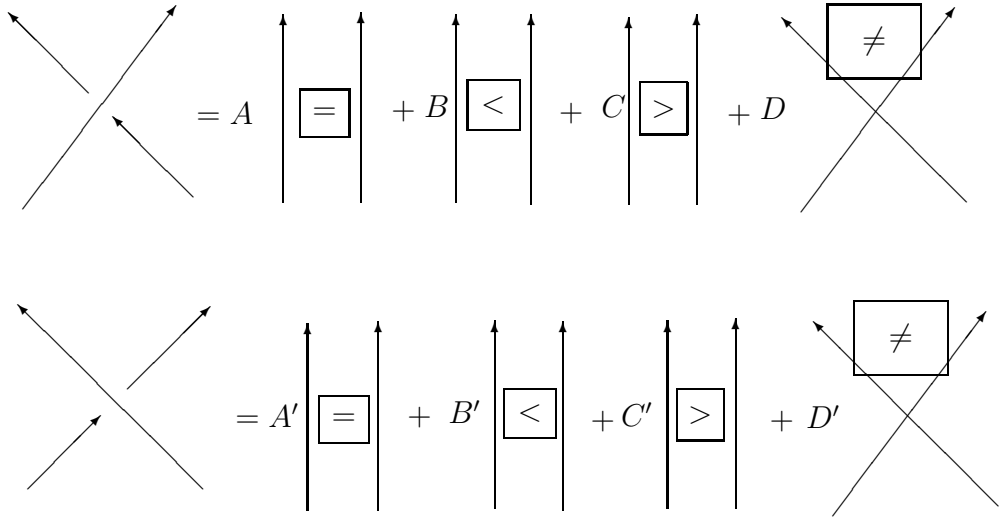

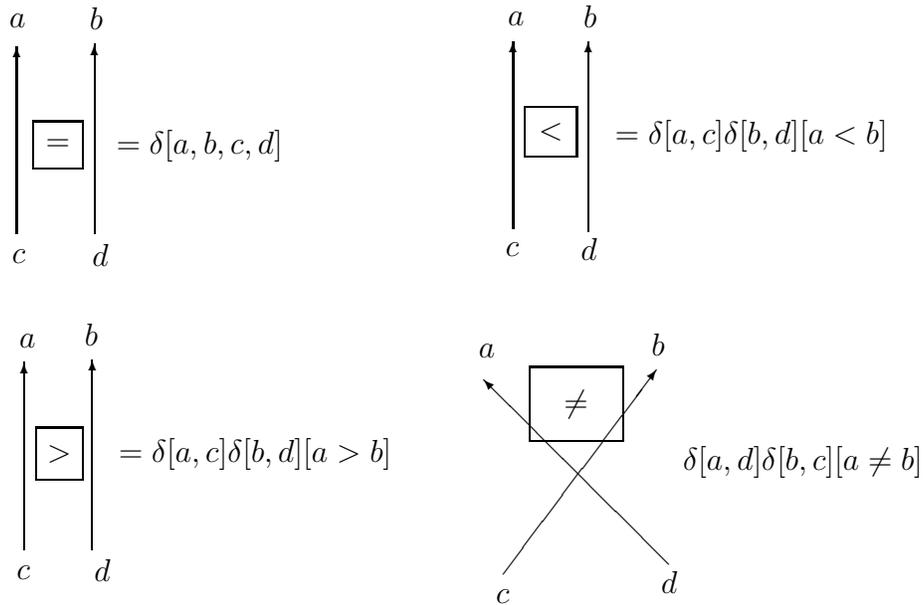
\begin{figure}
     \begin{center}
     \begin{tabular}{c}
{\tt    \setlength{\unitlength}{0.92pt}
\begin{picture}(381,269)
\thinlines    \put(224,80){\framebox(38,30){$\ne$}}
              \put(13,165){\vector(0,1){77}}
              \put(45,165){\vector(0,1){78}}
              \put(248,166){\vector(0,1){78}}
              \put(217,167){\vector(0,1){77}}
              \put(44,35){\vector(0,1){78}}
              \put(16,35){\vector(0,1){77}}
              \put(281,29){\vector(-1,1){76}}
              \put(213,25){\vector(3,4){63}}
              \put(222,197){\framebox(21,20){$<$}}
              \put(21,64){\framebox(19,21){$>$}}
              \put(20,192){\framebox(20,19){$=$}}
              \put(10,250){$a$}
              \put(43,249){$b$}
              \put(11,153){$c$}
              \put(44,152){$d$}
              \put(54,198){$= \delta[a,b,c,d]$}
              \put(215,251){$a$}
              \put(246,250){$b$}
              \put(214,156){$c$}
              \put(245,154){$d$}
              \put(259,202){$= \delta[a,c] \delta[b,d] [a < b]$}
              \put(14,118){$a$}
              \put(41,119){$b$}
              \put(13,23){$c$}
              \put(45,22){$d$}
              \put(55,71){$= \delta[a,c] \delta[b,d] [a > b]$}
              \put(203,114){$a$}
              \put(274,115){$b$}
              \put(210,13){$c$}
              \put(278,17){$d$}
              \put(287,68){$ \delta[a,d] \delta[b,c] [a \ne b]$}
\end{picture}}
     \end{tabular}
     \caption{\bf Local Matrices}
     \label{kroneckers}
\end{center}
\end{figure}

Now view Figure~\ref{homflypt} and Figure~\ref{kroneckers}. These figures give the conventions for setting up spin-preserving $R$--matrices from an index set $I.$
Each glyph in the first figure is defined as a combintation of Kronecker deltas in the second figure. The coefficents $A,,B,C,D,A',B',C',D'$ can be specified so that the indicated $R$--matrices satisfy the
Yang-Baxter equation. In particular, we can take $\lambda = q$ so that 
$$M_{ab} = q^{a/2} \delta_{ab}$$
and
$$A = q, B = (q-q^{-1}), C = 0, D=1$$ and $$A' = q^{-1}, B=0, C' = (q^{-1} - q), D'=1.$$ With this choice of coeffiecients we obtain, for classical links, a regular isotopy model of the Homflypt polynomial for each natural number $n= 1,2,3, \cdots,$ using the index set
$I_{n} = \{ -n, -n+2, \cdots , n-2, n \}.$  Letting $H_{n}[K](q)$ denote this model for a fixed value of $n,$ we have the following behaviour under classical curls (See \cite{KP,Turaev}) and crossings:
$$H_{n}[\Rcurl]=(q^{n+1}) H[\Arc],$$ 
$$H_{n}[\Lcurl] =(q^{-n-1}) H[\Arc],$$  
and the (regular isotopy) skein relation
$$H_{n}[\PosCross] - H_{n}[\NegCross] = (q-q^{-1})H_{n}[\OrSmooth].$$

By the same token, the model in the last paragraph extends by our procedure to a Homflypt polynomial for rotational virtuals. Lets call the regular isotopy version of this
invariant by the name $RH_{n}[K](q).$ This is the (natural) roational virtual extension of the Homflypt polynomial. More work needs to be done with this model. The reader will find that it
does not detect the difference (if there is one) between the links of Figure~\ref{linkcalc1} and Figure~\ref{linkcalc2} and corresponding unlinks. For such detection we shall perhaps need other quantum link invariants.\\

On the other hand there is a very important problem associated with extending the Homflypt polynomial to virtual knots. It is known \cite{Kho1,Kho2,Kho} that the Khovanov-Rozansky categorification of the invariants
$H_{n}(q)$ does extend to a well-defined homology theory for virtual knots and links. This means that the graded Euler characteristic of this homology theory gives an extension of the Homflypt polynomial to virtual knots and links. It is not known as of this writing whether such an extension can be achieved using models of the type discussed here with a modification of the virtual crossing operator (or some other elementary modification).\\

In the last paragraph we mention the possibility of modifiying the extensions of the Homflypt polynomial by modifiying the operator for the virtual crossing. This, in fact, is exactly what one can do in the case of the Alexander polynomial, obtaining a generalized Alexander polynomial in the context of quantum link invariants in \cite{KRO,KRCAT,GEN}. This invariant can be further modified to a non-trivial invariant
of rotational virtuals by using a standard virtual crossing operator, and we shall investigate that in a sequel to the present paper.\\

\section {Rotational Virtual Links, Quantum Algebras, Hopf algebras and the Tangle Category}

This section will show how the ideas and methods of this paper fit together with representations
of quantum algebras (to be defined below), Hopf algebras and invariants of virtual links.
We begin with a construction of the virtual tangle category. This category is a natural generalization of the virtual braid group. 
A functor from the virtual tangle category to an algebraic category will form a generalization of representations of virtual braid groups. This functor
is directly related to rotational invariants of virtual knots and links.  The category that we define contains virtual crossings, special elements that satisfy the algebraic Yang-Baxter equation and also cup and cap operators. The subcategory 
without the cup and cap operators and without any (symbolic) algebra elements except those involved with the algebraic Yang-Baxter operators is isomorphic to the String Category of \cite{CVBraid}.\\

A word to the reader about this section: In one sense this section is a review of known material in the form that Kauffman and Radford \cite{KRH} have shaped the theory of quantum invariants of knots and three-manifolds via finite-dimensional Hopf algebras. On the other hand, this theory is  generalized here to invariants of rotational virtual knots and links. This generalization is directly related to the structure of the virtual braid group as described in \cite{CVBraid}. More directly, the functor $F$ we describe  from the virtual tangle category to the category $Cat(A)$ of a Hopf algebra $A$, is quite natural, taking the virtual crossings to permutation operators in the Hopf algebra category. From the point of view of the Hopf algebra category, the virtual tangle category is its natural topological 
preimage.\\

Later in this section we shall show how this approach via categories and quantum algebras illuminates the structure of invariants that we have already described via state summations. In particular, we 
show how the link of Figure~\ref{linkcalc2} gets a non-trivial functorial image, corresponding to its non-trivial rotational bracket invariant. We also show that the link of Figure~\ref{linkcalc3} has a trivial functorial image. This means that not only is this link not detected by the rotational bracket polynomial, it is not detected by any quantum invariant formulated as outlined in this section. \\

\subsection{The Virtual Tangle Category}
The advantage in studying virtual knots up to rotational equivalence is that all quantum link
invariants generalize to invariants of rotational equivalence. This means that virtual rotational equivalence is
a natural equivalence relation for studying topology associated with solutions to the Yang-Baxter
equation. 
\smallbreak

\begin{figure}
     \begin{center}
     \begin{tabular}{c}
     \includegraphics[width=8cm]{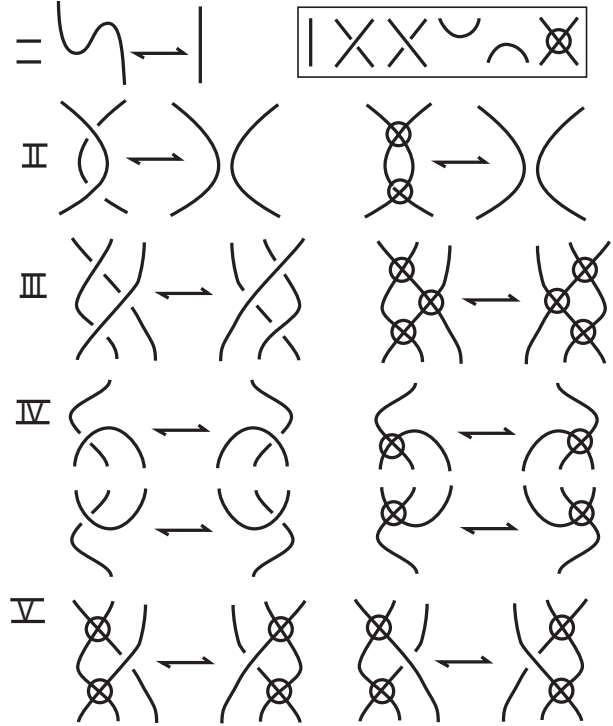}
     \end{tabular}
     \caption{Regular isotopy with respect to the vertical direction}
     \label{regtang}
\end{center}
\end{figure}

Here we create a context  by defining the {\it Virtual Tangle Category}, $VTC,$ as indicated
in Figure~\ref{regtang}. 
The tangle category is generated by the morphisms
shown in the box at the top of this figure. These generators are: a single identity line, right-handed and left-handed
crossings, a cap and a cup, a virtual crossing. The objects in the tangle category consist in the set of
$[n]$'s  where $n = 0, 1,2, \ldots.$ For a morphism $[n] \longrightarrow [m]$, the numbers $n$ and $m$ denote, respectively, the number of free arcs at the bottom and at the top of the diagram that represents  the morphism. The morphisms are like braids except that they can (due to the presence of the cups and caps) have different numbers of free ends at the top and the bottom of their diagrams. 
\smallbreak

\begin{figure}
     \begin{center}
     \begin{tabular}{c}
     \includegraphics[width=6cm]{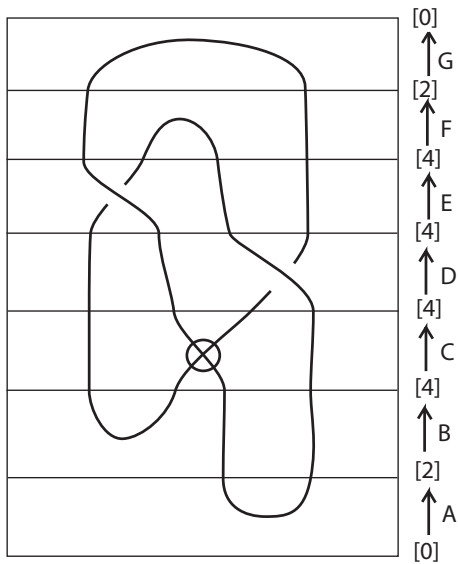}
     \end{tabular}
     \caption{Virtual trefoil as a morphism in the tangle category}
     \label{tangprod}
\end{center}
\end{figure}

The sense in which the elementary morphisms (line, cup, cap, crossings)  generate the tangle category is composition as  shown in Figure~\ref{tangprod}. For composition, the segments are matched so that the number of lower free ends on each segment is equal to the number of upper free ends on the segment below it.
The  Figure~\ref{tangprod} shows a virtual trefoil as a morphism from
$[0]$ to $[0]$ in the category. 
The tensor product of morphisms is the horizontal juxtaposition of their
diagrams. Each of the seven horizontal segments of the figure represents one of the elementary morphisms tensored  with the identity line. Consequently there is a well-defined composition of all of the segments and this composition is a morphism $[0] \longrightarrow [0]$ that represents the knot. 
\smallbreak

The basic equivalences of morphisms are
shown in Figure~\ref{regtang}. Note that in the type $IV$, the moves involving the minimum can be deduced from the other moves, and that we only indicate representative choices of crossing so that all crossings can be switched in any move to obtain another legal move. Note that $II, III, V$ are formally equivalent to the rules for unoriented
virtual braids. The zero-th move is a cancellation of consecutive maxima and minima, and the move $IV$
is a swing move in both virtual and classical relations of crossings to maxima and minima. It should be
clear that the tangle category is a generalization of the virtual braid group with a natural inclusion of unoriented virtual braids as special tangles in the category. Standard braid closure and the plat closure of braids
have natural definitions as tangle operations. Any virtual knot or link can be
represented in the tangle category as a morphism from $[0]$ to $[0],$ and one can prove that {\it two virtual links are
virtually regularly isotopic if and only if their tangle representatives are equivalent in the tangle category.}
None of the rules for equivalence in the tangle category involve either a classical loop or a
virtual loop. This means that the virtual tangle category is a natural home for the theory of rotational
virtual knots and links.
\smallbreak

\subsection{Quantum Algebra and Category}
Now we shift to a category associated with an  algebra that is directly related to our representations
of the virtual braid group. We take the following definition \cite{KP,KRH}: A {\it quantum algebra} $A$
is an algebra over a commutative ground ring $k$ with an invertible  mapping $s: A \longrightarrow A$
that is an {\it antipode}, that is $s(ab) = s(b)s(a)$ for all $a$ and $b$ in $A,$ and there is an element
$\rho \in A \otimes A$ satisfying the algebraic Yang-Baxter equation as in the equation below:
$$\rho_{12} \rho_{13} \rho_{23} = \rho_{23} \rho_{13} \rho_{12}.$$
We further assume that $\rho$ is invertible and that
$$ \rho^{-1} = (1_{A} \otimes s)\circ \rho = (s \otimes 1_{A})\circ \rho.$$
The multiplication in the algebra is usually denoted by $m: A\otimes A \longrightarrow A$ and is
assumed to be associative. It is also assumed that the algebra has a multiplicative unit element.
The defining properties of a quantum algebra are part of the properties of a Hopf algebra, but a Hopf algebra has a comultiplication $\Delta: A \longrightarrow A \otimes A$
that is a homomorphism of algebras, plus a list of further relations, including a fundamental relationship between the multiplication, the comultiplication and the antipode. In the interests of simplicity, we shall restrict ourselves to quantum algebras here, but most of the remarks that follow apply to Hopf algebras, and particularly quasi-triangular Hopf algebras. Information on Hopf algebras is included at the end of this section.  See \cite{KRH} for more about these connections.
\smallbreak

\begin{figure}
     \begin{center}
     \begin{tabular}{c}
     \includegraphics[width=8cm]{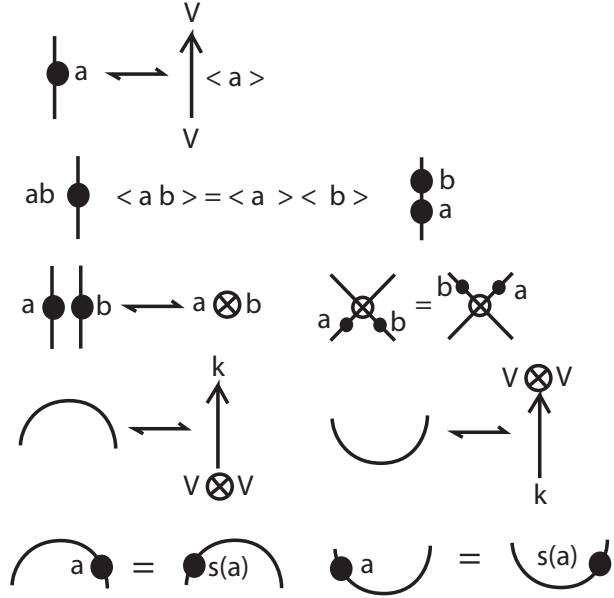}
     \end{tabular}
     \caption{ Morphisms in $Cat(A)$}
     \label{CatMorph}
\end{center}
\end{figure}

We construct a category $Cat(A)$ associated  with a quantum algebra $A$. This category is a very close relative to the virtual tangle category. $Cat(A)$ differs from the tangle category in that it has only virtual crossings, and
there are labeled vertical lines that carry elements of the algebra $A.$  See Figure~\ref{CatMorph}.
Each such labeled line is a morphism in the category.  The virtual crossing is a generating morphism as are the cups, caps and labeled lines. The objects in this category are the same
entities $[n]$ as in the tangle category. This category is identical in its framework to the tangle category but the crossings are not
present and lines labeled with algebra are present.
Given $a,b \in A$ we compose the morphisms corresponding to $a$ and $b$ by taking a line labeled
$ab$ to be their composition. In other words, if $\langle x \rangle$ denotes the morphism in
$Cat(A)$ associated with $x \in A$, then
$$\langle a \rangle \circ \langle b \rangle = \langle ab \rangle.$$
As for the additive structure in the algebra, we extend the category to an additive category by formally adding the generating morphisms (virtual crossings, cups, caps and algebra line segments). In Figure~\ref{CatMorph} we illustrate the composition of
such morphisms and we illustrate a number of other defining features of the category $Cat(A).$
\smallbreak

In the same figure we illustrate how the tensor product of elements $a \otimes b$ is represented by parallel vertical
lines with $a$ labeling the left line and $b$ labeling the right line. We indicate that the virtual
crossing acts as a permutation in relation to the tensor product of algebra morphisms. That is, we
illustrate that $$\langle a \rangle \otimes \langle b \rangle \circ P = P \circ \langle b \rangle \otimes \langle a \rangle.$$   Here $P$ denotes the virtual crossing of two segments, and is regarded as a morphism $P: V \otimes V \longrightarrow  V \otimes V$ (see remark below). Since the lines interchange, we expect $P$ to behave as the permutation of the two tensor factors.
\smallbreak

\begin{figure}
     \begin{center}
     \begin{tabular}{c}
     \includegraphics[width=8cm]{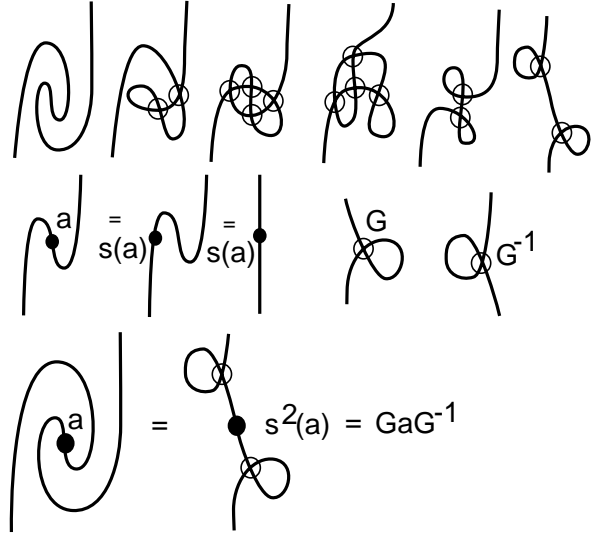}
     \end{tabular}
     \caption{Diagrammatics of the antipode}
     \label{antipode}
\end{center}
\end{figure}

In Figure~\ref{CatMorph} we show the notation
$V$ for the object $[1]$ in this category and we use $V \otimes V = [2]$, $V \otimes V \otimes V= [3]$
and so on for all the natural number objects in the category. We write
$[0] = k$, identifying the ground ring with the ``empty object" $[0].$ It is then axiomatic that
$k \otimes V = V \otimes k = V.$ Morphisms are indicated both diagrammatically and in terms of arrows
and objects in this figure. Finally, the figure indicates the arrow and object forms of the cup and the
cap, and crucial axioms relating the antipode with the cup and the cap.
A cap is regarded as a morphism from $V \otimes V$ to $k$, while
a cup is regarded as a morphism form $k$ to $V \otimes V.$ 
The basic property of the cup and the cap is the  {\em Antipode Property: if
one ``slides" a decoration across the maximum or minimum in a counterclockwise
turn, then the antipode $s$ of the algebra is applied to the decoration.}
In categorical terms this property says $$Cap \circ (\langle 1 \rangle \otimes a) = Cap \circ
(\langle sa \rangle \otimes 1 )$$ and $$(\langle a \rangle \otimes 1) \circ Cup  =  (1 \otimes \langle sa \rangle ) \circ Cup.$$  Here $1$ denotes the identity morphism for $[0]$. These properties and other naturality properties of the cups and the
caps are illustrated in Figure~\ref{CatMorph} and Figure~\ref{antipode}.  The naturality properties of the flat diagrams  in this category include regular homotopy of immersions (for diagrams without algebra decorations), as illustrated in these figures.
\smallbreak

In Figure~\ref{antipode}  we see how the antipode property of the cups and caps leads to a
diagrammatic interpretation of the antipode. In the figure we see that the antipode $s(a)$ is represented by composing with a cap and a cup on either side of the morphism for $a$. In terms of the composition
of morphisms this diagram becomes
$$\langle sa \rangle = (Cap \otimes 1)  \circ (1 \otimes \langle a \rangle \otimes 1)\circ(1 \otimes  Cup).$$
Similarly, we have
$$\langle s^{-1}a \rangle = (1 \otimes Cap)  \circ(1 \otimes \langle a \rangle \otimes 1)\circ( Cup \otimes 1).$$
This, in turn, leads to the
interpretation of the flat curl as an  element $G$ in $A$ such that
$s^{2}(a) = GaG^{-1}$ for all $a$ in $A.$  $G$ is a flat curl diagram interpreted
as a morphism in the category. We see that, formally, it is natural to interpret
$G$ as an element of $A$.   In a so-called {\em ribbon Hopf
algebra} there is such an element  already in the algebra. In the general case it
is natural to extend the algebra to contain such an element.
\smallbreak

\begin{figure}
     \begin{center}
     \begin{tabular}{c}
     \includegraphics[width=8cm]{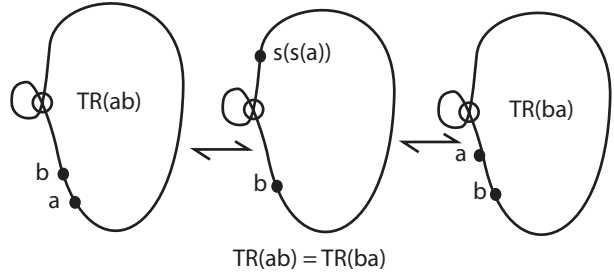}
     \end{tabular}
     \caption{Formal trace}
     \label{trace}
\end{center}
\end{figure}

\subsection{The Basic Functor and the Rotational Trace}
\begin{figure}
     \begin{center}
     \begin{tabular}{c}
     \includegraphics[width=8cm]{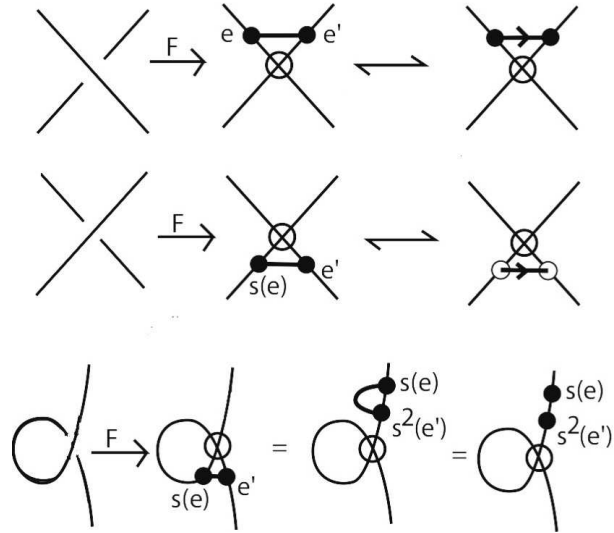}
     \end{tabular}
     \caption{The functor $F: VTC \longrightarrow Cat(A)$}
     \label{Functor}
\end{center}
\end{figure}

We are now in a position to describe a functor $F$ from the  virtual tangle category $VTC$
to $Cat(A).$  (Recall that the virtual tangle category is defined for virtual link diagrams without
decorations. It has the same objects as $Cat(A).$) 
$$F:VTC \longrightarrow Cat(A)$$
The functor $F$ decorates each
positive  crossing of the tangle (with respect to the vertical - see Figure~\ref{Functor})
with the Yang-Baxter element (given by the quantum algebra $A$)
$\rho = \Sigma e \otimes e^{'}$ and each negative crossing (with respect to the
vertical) with $\rho^{-1} = \Sigma s(e) \otimes e^{'}$. The form of the
decoration is indicated in Figure~\ref{Functor}. Since we have labelled the negative crossing with the
inverse Yang-Baxter element, it follows at once that the two crossings are mapped to inverse elements in the category of the algebra.  This association is a direct generalization of our mapping
of the virtual braid group to the stringy connector presentation \cite{CVBraid}.
\smallbreak

We now point out the structure of the image of a knot, link or tangle under this functor.
The key point about this functor is that, because quantum algebra elements can be
moved around the diagram, we can concentrate all the image algebra in one place.
Because the flat curls are identified with either $G$ or $G^{-1}$, we can use
regular homotopy of  immersions to bring the image under $F$ of each component of a virtual link diagram to the form of a circle with a single concentrated decoration (involving a sum over many
products) and a reduced pattern of flat curls that can be encoded as a power of the special
element $G.$ Once the underlying curve of a link  component is converted to a loop with
total turn zero, as in Figure~\ref{trace}, then we can think of such a loop, with algebra labeling the
loop, as a representative for a formal trace of that algebra and call it $TR(X)$ as in the figure.
In the figure we illustrate that for such a labeling $$TR(ab) = TR(ba),$$ thus one can take a product of
algebra elements on a zero-rotation loop up to cyclic order of the product. In situations where we
choose a representation of the algebra or in the case of finite dimensional Hopf algebras where one
can use right integrals \cite{KRH}, there are ways to make actual evaluations of such traces. Here we
use them formally to indicate the result of concentrating the algebra on the loop.
\smallbreak

\begin{figure}
     \begin{center}
     \begin{tabular}{c}
     \includegraphics[width=8cm]{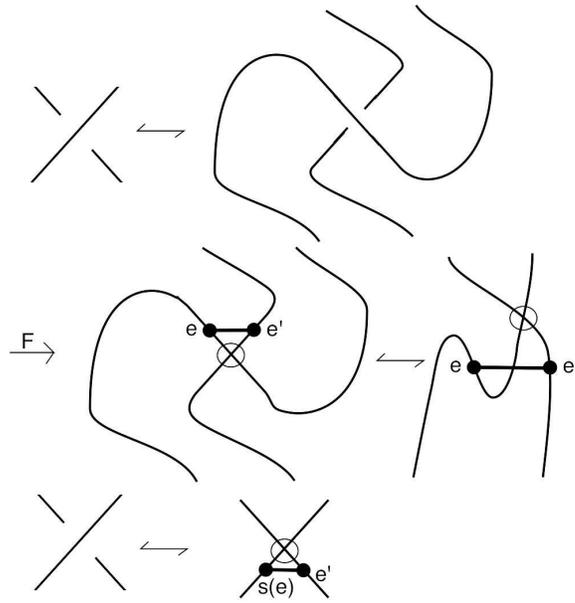}
     \end{tabular}
     \caption{Inverse and antipode}
     \label{twist}
\end{center}
\end{figure}

\begin{figure}
     \begin{center}
     \begin{tabular}{c}
     \includegraphics[width=9cm]{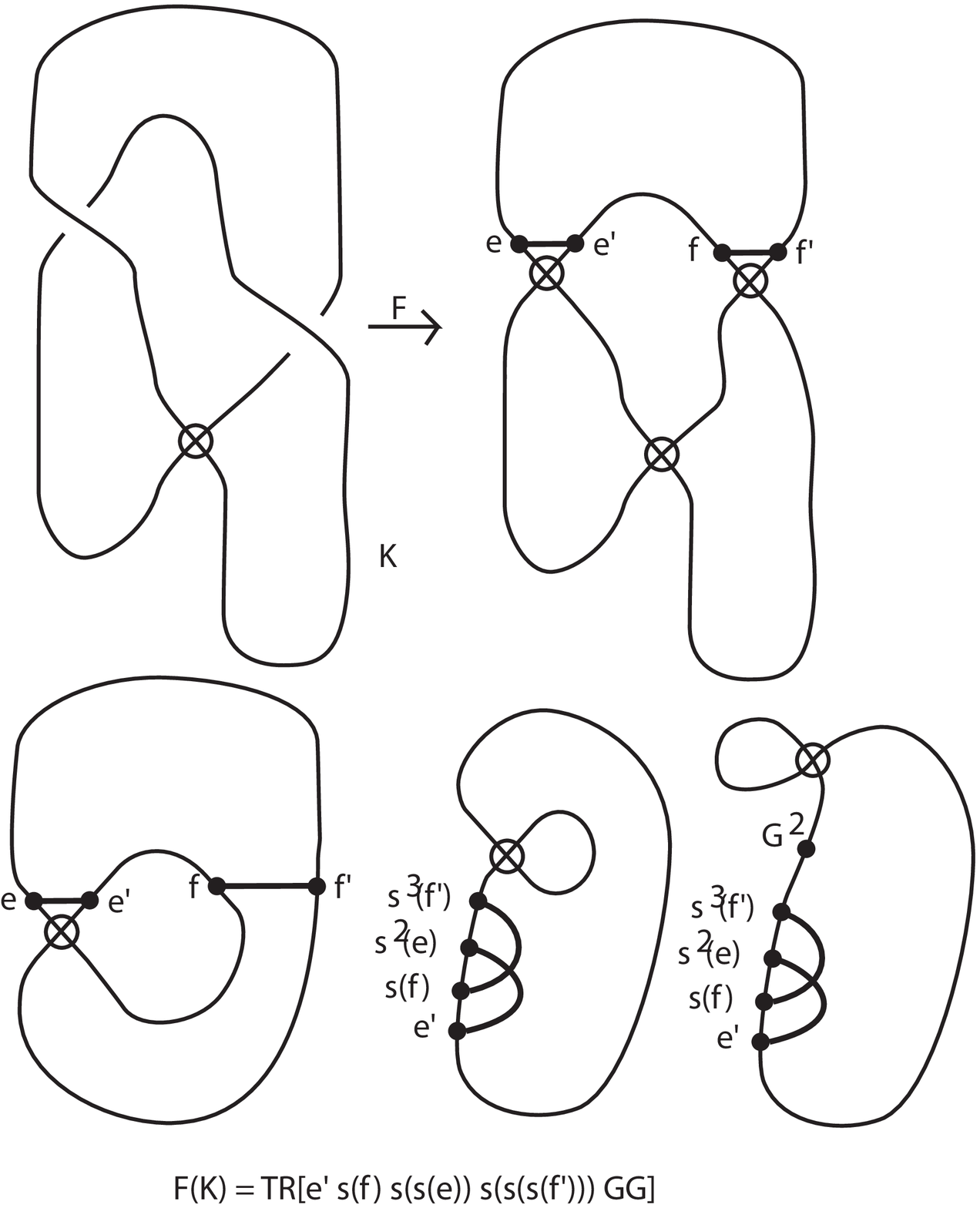}
     \end{tabular}
     \caption{The functor $F: T \longrightarrow Cat(A)$ applied to a virtual trefoil}
     \label{ApplyFunctor}
\end{center}
\end{figure}

One further comment is in order about the antipode. In Figure~\ref{twist} we show that our axiomatic assumption about the antipode (the sliding rule around maxima and minima) actually demands that
the inverse of $\rho$ is $(s \otimes 1_{A})\circ \rho = (1_{A} \otimes s )\circ \rho$. This follows by
examining the form of the inverse of the positive crossing in the tangle category by turning that crossing
to produce an identity between the positive crossing and the negative crossing twisted with additional
maxima and minima. This relationship shows that if we set the functor $F$ on a right-handed crossing as we have done, then the way it maps the inverse crossing is forced and that this inverse corresponds to the inverse of $\rho$ in the quantum algebra. Thus the quantum algebra formula for the inverse of 
$\rho$ is forced by the topology.
\smallbreak

In Figure~\ref{ApplyFunctor} we illustrate the entire functorial process for  the virtual trefoil of Figure~\ref{tangprod}. The virtual trefoil is denoted by $K$, and we find that $F(K)$ reduces to a zero-rotation circle with the
inscription $ e' s(f) s^{2}(e) s^{3}(f') G^{2}. $ We can, therefore, write the equation
$$F(K) = TR[e' s(f) s^{2}(e) s^{3}(f') G^{2}].$$
Another way to think about this trace expression is to regard it as a Gauss code for the knot that has extra structure.  The powers of the antipode and the power of $G$ keep track of rotational features in the diagram as it lives in the tangle category up to regular isotopy.
We now see that the mapping of the virtual braid group to the braid group generated by permutations
and string connectors has been generalized to the functor $F$ taking the virtual tangle category
to the abstract category of a quantum algebra. We regard this generalization as an appropriate context for thinking about virtual knots, links and braids.
\smallbreak

The category $Cat(A)$ of a quantum algebra $A$ can be generalized to an abstract category
with labels, virtual crossings, and with stringy connections that satisfy the algebraic Yang-Baxter equation. Each such stringy connection has a left label $e$ or $s(e)$ and a right label $e'.$ We retain the formalism of the antipode
as a formal replacement for adjoining a label with a cup and a cap. The resulting {\it abstract algebra category} will be denoted by $\overline{Cat(A)}$. Since we take this category with no further relations,
the functor $ \overline{F}:VTC \longrightarrow \overline{Cat(A)}$ is an equivalence of categories. This functor is
the direct analog of our reformulation of the virtual braid group in terms of stringy connectors in \cite{CVBraid}.

\subsection{Applying The Functor to Rotational Virtual Links}
We now show how this approach via categories and quantum algebras illuminates the structure of invariants that we have already described via state summations. In particular, we 
show how the link of Figure~\ref{linkcalc2} gets a non-trivial functorial image, corresponding to its non-trivial rotational bracket invariant. We also show that the link of Figure~\ref{linkcalc3} has a trivial functorial image. This means that not only is this link not detected by the rotational bracket polynomial, it is not detected by any quantum invariant formulated as outlined in this section.\\

\begin{figure}
     \begin{center}
     \begin{tabular}{c}
     \includegraphics[width=8cm]{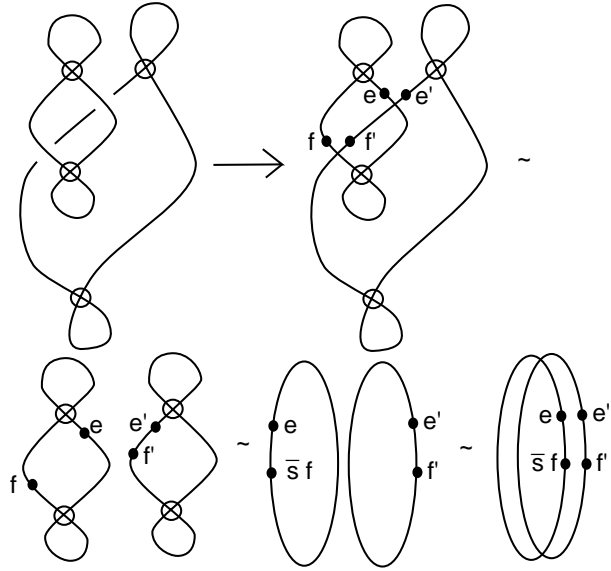}
     \end{tabular}
     \caption{First Rotational Link - Non-Trivial Invariants}
     \label{nontriv}
\end{center}
\end{figure}

In Figure~\ref{nontriv} we apply the functor (to the category of abstract quantum algebras) to the link $L_{2}$ of Figure~\ref{linkcalc2}. We have already shown that this link is detected by the rotational bracket. We see from Figure~\ref{nontriv} that there is no apparent reduction of the diagram in abstract quantum algebra category. The juxtaposed algebra elements are not inverses of one another.
Recall that if $\rho$ denotes the algebraic Yang-Baxter element, then $\rho^{-1} = (s \otimes 1)\rho = (1 \otimes s^{-1})\rho.$ Such inverses do not occur in this diagram.\\

\begin{figure}
     \begin{center}
     \begin{tabular}{c}
     \includegraphics[width=8cm]{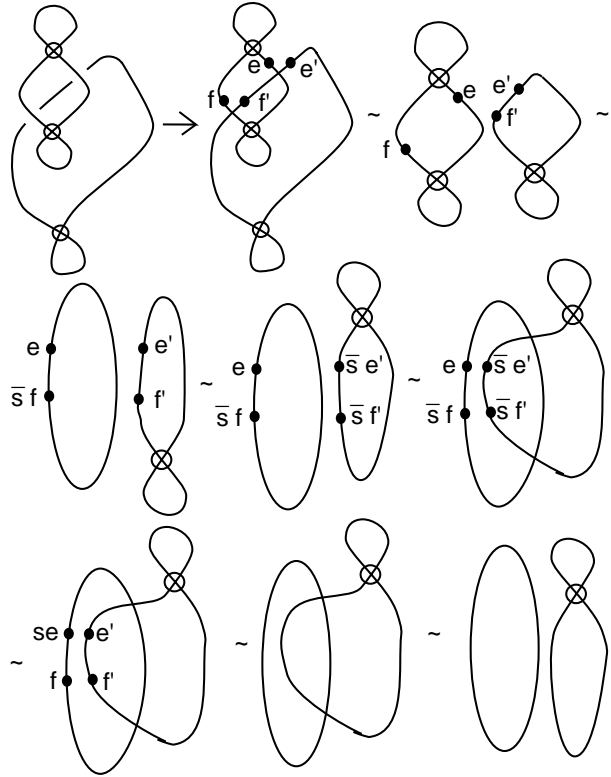}
     \end{tabular}
     \caption{Second Rotational Link - Trivial Invariants}
     \label{triv}
\end{center}
\end{figure}

On the other hand, in Figure~\ref{triv} we apply the functor to the link $L_{3}$ of Figure~\ref{linkcalc3}, and we find that by suitable transfomations of the diagram, the functor will evaluate this link
just as though it were an unlink of the corresponnding components of $L_{3}.$ {\it Thus no quantum link invarinant of the type considered in this section can detect $L_{3}.$} \\

The next two subsections give extra information about Kirby calculus and Hopf algebras. We can, study rotational links up to the equivalence relation generated by the Kirby calculus. This is a project for another paper.\\

\subsection{Hopf Algebras and Kirby Calculus}
\begin{figure}
     \begin{center}
     \begin{tabular}{c}
     \includegraphics[width=8cm]{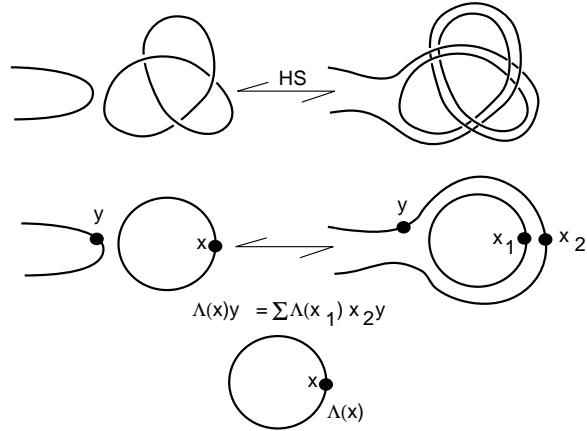}
     \end{tabular}
     \caption{The Kirby move}
     \label{kirby}
\end{center}
\end{figure}

In Figure~\ref{kirby} we illustrate how one can use this concentration of algebra on the loop in the context of a Hopf algebra that has a right integral. The right integral is a function
$\lambda: A \longrightarrow k$ satisfying $$\lambda(x) 1_{A} = \Sigma \lambda(x_{1})x_{2}$$
where the coproduct in the Hopf algebra has formula $\Delta(x) = \Sigma x_{1} \otimes x_{2}$.
Here we point out how the use of the coproduct corresponds to doubling the lines in the diagram, and that if one were to associate the function
$\lambda$ with a circle with rotation number one,  then the resulting link evaluation will be invariant under the so-called Kirby move \cite{KCalc}.
The Kirby move replaces two link components with new ones by doubling one component and connecting one of the components of the double with the other component. Under our functor from the virtual tangle category to the category for the Hopf algebra, a knot goes to a circle with algebra concentrated at $x.$ The doubling of the knot goes to concentric circles labeled with the coproduct 
$\Delta(x) = \Sigma x_{1} \otimes x_{2}.$  Figure~\ref{kirby} shows how invariance under the handle-slide in the tangle category corresponds the integral equation
$$\lambda(x) y = \Sigma \lambda(x_{1})x_{2} y.$$
It turns out that
classical framed links $L$ have an associated compact oriented three manifold $M(L)$ and that
two links related by Kirby moves have homeomorphic three-manifolds. Thus the evaluation of links using the right integral yields invariants of three-manifolds. Generalizations to virtual three-manifolds are under investigation \cite{DK1}. We only sketch this point of view here, and refer the reader to \cite{KRH}.
\smallbreak

\subsection {Hopf Algebra}

This section is added for reference about Hopf algebras. Quasitriangular Hopf algebras are an
important special case of the quantum algebras discussed in this section.
\smallbreak

Recall that a {\it Hopf algebra} \cite{Sweedler} is a bialgebra  $A$ over a commutative ring $k$ that has an associative multiplication $m:A \otimes A \longrightarrow A,$
and a coassociative comultiplication, and is equipped with a counit, a unit
and an antipode. The ring $k$ is usually taken to be a field.  The associative law for the multiplication  $m$ is expressed by the equation $$m(m \otimes1_{A}) = m(1_{A} \otimes m)$$ where $1_{A}$ denotes the identity map on A.
\smallbreak

The coproduct
$\Delta :A \longrightarrow A \otimes A$ is an algebra homomorphism and is
coassociative in the sense that
$$(\Delta \otimes 1_{A})\Delta = (1_{A} \otimes \Delta) \Delta.$$ 

The {\it unit} is a mapping from $k$ to $A$ taking $1_k$ in $k$ to $1_A$ in $A$ and, thereby, defining an action of $k$ on $A.$  
It will be convenient to just identify the $1_k$ in $k$ and the $1_A$ in $A$, and to ignore
the name of the map that gives the unit. 
\smallbreak

The counit is an algebra mapping from $A$ to $k$ denoted by $\epsilon :A \longrightarrow k.$ The following formula for the counit dualize the structure inherent in the unit:
$$(\epsilon \otimes 1_{A}) \Delta = 1_{A} = (1_{A} \otimes \epsilon) \Delta.$$

It is convenient to write formally
$$\Delta (x) = \sum x_{1} \otimes x_{2} \in A \otimes A$$
to indicate the decomposition of the coproduct of $x$ into a sum of first and second factors in the two-fold tensor product of $A$ with itself. We shall often drop the summation sign and write
$$\Delta (x) = x_{1} \otimes x_{2}.$$

The antipode is a mapping $s:A \longrightarrow A$ satisfying the equations
$$m(1_{A} \otimes s) \Delta (x) = \epsilon (x)1_{A} \quad \mbox{and} \quad m(s \otimes 1_{A}) \Delta (x)= \epsilon (x)1_{A}.$$   
It is a consequence of this definition that $s(xy) = s(y)s(x)$ for all $x$ and $y$ in A. \vspace{3mm}

A {\it quasitriangular Hopf algebra} \cite{Drinfeld} is a Hopf algebra  $A$ with an element 
$\rho \in A \otimes A$ satisfying
the following conditions:
\vspace{3mm}

\noindent
(1) $\rho \Delta = \Delta' \rho$ where $\Delta'$ is the composition
of $\Delta$ with the map on
$A \otimes A$ that switches the two factors. \vspace{3mm}

\noindent
(2) $$\rho_{13} \rho_{12} = (1_{A} \otimes \Delta) \rho,$$ $$\rho_{13} \rho_{23} = (\Delta \otimes 1_{A})\rho.$$

The symbol $\rho_{ij}$ denotes the placement of the first and second tensor factors of $\rho$ in the $i$ and $j$ places in a triple tensor product. For example, if $\rho = \sum e \otimes e'$ then $$\rho_{13} = \sum e \otimes 1_{A} \otimes e'.$$

Conditions (1) and (2) above imply that $\rho$ has an inverse and that

$$ \rho^{-1} = (1_{A} \otimes s^{-1}) \rho = (s \otimes 1_{A}) \rho.$$

It follows easily from the axioms of the quasitriangular Hopf algebra that $\rho$ satisfies the Yang-Baxter equation

$$\rho_{12} \rho_{13} \rho_{23} = \rho_{23} \rho_{13} \rho_{12}.$$

A less obvious fact about quasitriangular Hopf algebras is that there exists an element $u$ such that $u$ is invertible and $s^{2}(x) = uxu^{-1}$ for all $x$ in $A.$ In fact, we may take $u = \sum s(e')e$ where $\rho = \sum e \otimes e'.$ This result, originally due to Drinfeld \cite{Drinfeld}, follows from the diagrammatic categorical context of this paper. \vspace{3mm}

An element $G$ in a Hopf algebra is said to be {\em grouplike} if $\Delta (G) = G \otimes G$ and $\epsilon (G)=1$ (from which it follows that $G$ is invertible and $s(G) = G^{-1}$). A quasitriangular Hopf algebra is said to be a {\em ribbon Hopf algebra} \cite{RTG,KRH} if there exists a grouplike element $G$ such that (with $u$ as in the previous paragraph) $v = G^{-1}u$ is in
the center of $A$ and $s(u) = G^{-1}uG^{-1}$. We call $G$ a special grouplike element of $A.$
\vspace{3mm}

Since $v=G^{-1}u$ is central, $vx=xv$ for all $x$ in $A.$ Therefore $G^{-1}ux = xG^{-1}u.$ We know that $s^{2}(x) = uxu^{-1}.$ Thus $s^{2}(x) =GxG^{-1}$ for all $x$ in $A.$ Similarly, $s(v) = s(G^{-1}u) = s(u)s(G^{-1})=G^{-1}uG^{-1}G =G^{-1}u=v.$ Thus, the square of the
antipode is represented as conjugation by the special grouplike element in a ribbon Hopf algebra, and the central element $v=G^{-1}u$ is invariant under the antipode.
\smallbreak

This completes the summary of Hopf algebra properties that are relevant to the last section of the paper.




\begin{thebibliography}{AF}


\bibitem{CS1}
J.S.Carter, S. Kamada and M. Saito, Stable equivalence of knots on surfaces and virtual  knot
cobordisms, in ``Knots 2000 Korea, Vol. 1 (Yongpyong)", {\em JKTR} {\bf 11}, No. 3 (2002), 311--320.

\bibitem{Drinfeld}
V.G. Drinfeld, Quantum groups,
{\em Proceedings of the International
Congress of Mathematicians, Berkeley, California, USA}(1987), 798--820. Amer. Math. Soc., Providence, RI.

\bibitem{Dye}
H. Dye, Characterizing Virtual Knots, Ph.D. Thesis (2002), UIC.

\bibitem{DK1} H. Dye and L. H. Kauffman, Virtual knot diagrams and the
Witten-Reshetikhin-Turaev invariant. math.GT/0407407, {\it J. Knot Theory Ramifications} 14 (2005), no. 8, 1045--1075.

\bibitem{DKK}
H. Dye, A.  Kaestner and L. H. Kauffman, Khovanov homology, Lee homology and a Rasmussen invariant for virtual knots.
arXiv:1409.5088.

\bibitem{GPV}
 M. Goussarov, M.Polyak and O. Viro, Finite type invariants of classical and virtual knots,
{\em Topology}  {\bf 39} (2000),  1045--1068.

 \bibitem{JO} 
V.F.R. Jones, A polynomial invariant for links via von Neumann algebras,
Bull. Amer. Math. Soc. {\bf 129} (1985), 103--112.

\bibitem{JO1}
 V.F.R.Jones.  Hecke algebra representations of braid groups and link polynomials.  Ann. of Math.  126 (1987), pp. 335-338.

\bibitem{JO2}
V.F.R.Jones.  On knot invariants related to some statistical mechanics models.  Pacific J. Math., vol. 137, no. 2 (1989), pp. 311-334.

\bibitem{KaestnerKauffman} A. M. Kaestner.; L. H. Kauffman, Parity, skein polynomials and categorification. J. Knot Theory Ramifications 21 (2012), no. 13, 1240011, 56 pp.

\bibitem{KaB}  L.H. Kauffman, (1987), State Models and the Jones Polynomial,
{\bf Topology} {\bf 26}, pp. 395-407.

\bibitem{VKT} L. H. Kauffman, Virtual Knot Theory , {\em European J. Comb.} (1999) Vol. 20, 663-690.

\bibitem{SVKT}
L. H. Kauffman, A Survey of Virtual Knot Theory in {\em Proceedings of Knots in Hellas '98},
World Sci. Pub. 2000 , pp. 143-202.

\bibitem{DVK}
L.H. Kauffman, Detecting Virtual Knots, in Atti. Sem. Mat. Fis.
Univ. Modena Supplemento al Vol. IL, 241-282 (2001).

\bibitem{SL} L. H. Kauffman, A self-linking invariant of virtual
knots. {\it Fund. Math.} 184 (2004), 135--158.

\bibitem{CVBraid}L. H. Kauffman; S. Lambropoulou, A categorical model for the virtual braid group. J. Knot Theory Ramifications 21 (2012), no. 13, 1240008, 48 pp.
 
\bibitem{VKCob} L. H. Kauffman, Virtual knot cobordism, in ``New Ideas in Low Dimensional Topology", edited by L. H. Kauffman and V. O. Manturov, Knots and Everything Series, Vol. 56, World Scientific Pub. Co. (2014), pp. 335-377.

\bibitem {KP}
L.H. Kauffman, {\em Knots and Physics}, World Scientific Publishers (1991), 
Second Edition (1993), Third Edition (2002), Fourth Edition (2012).

\bibitem{KRH}
L.H. Kauffman and D.E. Radford, Invariants of 3- manifolds derived from finite dimensional Hopf algebras, {\em JKTR} {\bf 4}, No. 1 (1995), 131--162.

\bibitem{FKT} L. H. Kauffman, {\em Formal Knot Theory}, Princeton
University Press, Lecture Notes Series  30 (1983). Dover Publications (2005).

\bibitem{KRO}
L. H. Kauffman and David E. Radford, Oriented quantum algebras and
invariants of knots and links,  {\em Journal of Algebra}, Vol. 246, 253-291 (2001).

\bibitem{KRCAT}
L. H. Kauffman and David E. Radford, Oriented quantum algebras,
categories and invariants of knots and links.  {\em JKTR}, vol 10, No. 7 (2001), 1047-1084. 

\bibitem{GEN}
L.H. Kauffman and D.E. Radford, Bi-oriented Quantum Algebras, and a
Generalized Alexander Polynomial for Virtual Links, in  
``Diagrammatic Morphisms and Applications" (San Francisco, CA, 2000), 113--140, Contemp. Math., 318, Amer. Math. Soc., Providence, RI,
2003.
   
\bibitem{BraidGates}
L.H. Kauffman and Samuel J. Lomonaco, Braiding Operators are Universal Quantum
Gates, New Journal of Physics 6 (2004) 134, pp. 1-39.


\bibitem{KCalc} R. Kirby, A calculus for framed links in $S^{3}$, {\em Invent.
Math.}, {\bf 45} (1978), 35--56.

\bibitem{LOMI} N.Y. Reshetikhin, Quantized universal enveloping algebras, the
Yang-Baxter equation and invariants of links, I and II, {\em  LOMI reprints
E-4-87 and E-17-87}, Steklov Institute, Leningrad, USSR.

\bibitem{RTG}
N. Yu. Reshetikhin and V.G. Turaev,
Ribbon graphs and their invariants derived from quantum groups.
{\em Commun. Math. Phys. }  {\bf 127} (1990), 1--26.

\bibitem{Turaev} V.G. Turaev, The Yang-Baxter equations and invariants of links,
LOMI preprint E-3-87, Steklov Institute, Leningrad, USSR., {\em Inventiones
Math.},Vol. 92, Fasc.3, pp. 527-553.


\bibitem{Sweedler} M.E. Sweedler, ``Hopf Algebras", Mathematics Lecture Notes Series,
Benjamin, New York, 1969.


\bibitem{Hennings} M.A. Hennings, Hopf algebras and regular isotopy invariants
for link diagrams, {\em Math. Proc. Camb. Phil. Soc.}, Vol. 109 (1991), pp. 59-77


\bibitem{Kho2}
 M. Khovanov ; L. Rozansky, Matrix factorizations and link homology. II. Geom. Topol. 12 (2008), no. 3, 1387Ð1425.
 
 \bibitem{Kho1}
 M. Khovanov ; L. Rozansky, Matrix factorizations and link homology. Fund. Math. 199 (2008), no. 1, 1Ð91. 
 
 \bibitem{Kho}
 M. Khovanov ; L. Rozansky, Virtual crossings, convolutions and a categorification of the SO(2N) Kauffman polynomial. 
 

\bibitem{Ma}
V.O. Manturov and D. P. Ilyutko , ``Virtual Knots -- The State of the Art," Knots and Everything Series, Vol. 51, World Scientific Pub. Co. (2012).

\bibitem{MP1} V. O. Manturov, Parity and Cobordisms of Free Knots ( to appear in Sbornik),
arXiv:1001.2827 

\bibitem{PCFree} V. O. Manturov,  Parity and cobordisms of free knots. (Russian) Mat. Sb. 203 (2012), no. 2, 45--76; translation in Sb. Math. 203 (2012), no. 1-2, 196Ð223.

\bibitem{CobFree} D. P. Ilyutko; V. O. Manturov, Cobordisms of free knots. (Russian) Dokl. Akad. Nauk 429 (2009), no. 4, 439--441; translation in Dokl. Math. 80 (2009), no. 3, 844Ð846.


\bibitem{MP} V. O.~Manturov, Parity in knot theory', {\it Math.\ sb.} {\bf
201}:5, 693--733 (2010) (Original Russian Text in {\it Mathematical
sbornik} {\bf 201}:5, pp. 65--110 (2010)).

\bibitem{Witten} E. Witten. Quantum Field Theory and the Jones Polynomial. Comm. in Math. Phys.
Vol. 121 (1989), 351-399.


\end{thebibliography}
\end{document}